# CHARACTERIZATION ALGORITHMS FOR SHIFT RADIX SYSTEMS WITH FINITENESS PROPERTY

MARIO WEITZER

*Chair of Mathematics and Statistics, Montanuniversität Leoben, Franz Josef-Straße 18*
*Leoben, 8700, Austria*
*mario.weitzer@unileoben.ac.at*

For $d \in \mathbb{N}$ and $\mathbf{r} \in \mathbb{R}^d$ let $\tau_{\mathbf{r}} : \mathbb{Z}^d \to \mathbb{Z}^d$, where $\tau_{\mathbf{r}}(\mathbf{a}) = (a_2, \ldots, a_d, -\lfloor \mathbf{ra} \rfloor)$ for $\mathbf{a} = (a_1, \ldots, a_d)$, denote the *(d-dimensional) shift radix system associated with* $\mathbf{r}$. $\tau_{\mathbf{r}}$ is said to have the *finiteness property* iff all orbits of $\tau_{\mathbf{r}}$ end up in $(0, \ldots, 0)$; the set of all corresponding $\mathbf{r} \in \mathbb{R}^d$ is denoted by $\mathcal{D}_d^{(0)}$, whereas $\mathcal{D}_d$ consists of those $\mathbf{r} \in \mathbb{R}^d$ for which all orbits are eventually periodic. $\mathcal{D}_d^{(0)}$ has a very complicated structure even for $d = 2$.

In the present paper two algorithms are presented which allow the characterization of the intersection of $\mathcal{D}_d^{(0)}$ and any closed convex hull of finitely many interior points of $\mathcal{D}_d$ which is completely contained in the interior of $\mathcal{D}_d$. One of the algorithms is used to determine the structure of $\mathcal{D}_2^{(0)}$ in a region considerably larger than previously possible, and to settle two questions on its topology: It is shown that $\mathcal{D}_2^{(0)}$ is disconnected and that the largest connected component has nontrivial fundamental group. The other algorithm is the first characterizing $\mathcal{D}_d^{(0)}$ in a given convex polyhedron which terminates for all inputs. Furthermore several infinite families of “cutout polygons” are deduced settling the finiteness property for a chain of regions touching the boundary of $\mathcal{D}_2$.

## 0. Conventions

- $\mathbb{N}$ is defined as the set of all positive integers, and $\mathbb{N}_0 := \mathbb{N} \cup \{0\}$.
- The modulo function mod has precedence over addition and subtraction.
- $\mathcal{P}(M)$ denotes the power set of a set $M$.
- For $x \in \mathbb{R}$, $\lfloor x \rfloor$ (floor) denotes the largest integer less than or equal to $x$ and $\lceil x \rceil$ (ceiling) denotes the smallest integer greater than or equal to $x$. The floor or ceiling of a vector $\mathbf{r} \in \mathbb{R}^d$ is obtained by elementwise application.
- $\mathbf{0}$ means the zero vector of suitable length.
- $\mathbf{rs}$ is the inner product of $\mathbf{r} \in \mathbb{R}^d$ and $\mathbf{s} \in \mathbb{R}^d$.
- $\mathrm{id}_M$ is the identity map on a set $M$.

- For a subset $M$ of a topological space $\operatorname{int}(M)$ denotes the interior, $\overline{M}$ the closure, and $\partial M$ the boundary of M.

## 1. Introduction

In 2005 Akiyama *et al.* [1] introduced so-called *shift radix systems* (compare also [2,3,4]). For a natural number $d$ and a real vector $\mathbf{r} \in \mathbb{R}^d$ the mapping $\tau_{\mathbf{r}} : \mathbb{Z}^d \to \mathbb{Z}^d$ defined by

$$\tau_{\mathbf{r}}(\mathbf{a}) = (a_2, \ldots, a_d, -\lfloor \mathbf{ra} \rfloor) \qquad (\mathbf{a} = (a_1, \ldots, a_d)), \tag{1.1}$$

is called the *d-dimensional shift radix system associated with* $\mathbf{r}$ (*SRS*) and $\mathbf{r}$ its *parameter*. Let

$$\mathcal{D}_d := \left\{\mathbf{r} \in \mathbb{R}^d \mid \forall\, \mathbf{a} \in \mathbb{Z}^d : \exists\, (m,n) \in \mathbb{N}^2 : m \neq n \wedge \tau_{\mathbf{r}}^m(\mathbf{a}) = \tau_{\mathbf{r}}^n(\mathbf{a})\right\} \tag{1.2}$$

$$\mathcal{D}_d^{(0)} := \left\{\mathbf{r} \in \mathbb{R}^d \mid \forall\, \mathbf{a} \in \mathbb{Z}^d : \exists\, n \in \mathbb{N} : \tau_{\mathbf{r}}^n(\mathbf{a}) = \mathbf{0}\right\} \tag{1.3}$$

where for any $n \in \mathbb{N}_0$, $\tau_{\mathbf{r}}^n(\mathbf{a})$ means the $n$-fold iterative application of $\tau_{\mathbf{r}}$ to $\mathbf{a}$. The SRS $\tau_{\mathbf{r}}$ is said[a b] to have the *finiteness property* iff $\mathbf{r} \in \mathcal{D}_d^{(0)}$.

It is well-known that $\mathcal{D}_d^{(0)}$ has a very complicated structure even for $d = 2$ (cf. [18] or Fig. 2 in Section 4). In the present paper two algorithms are presented (Section 3) which allow the characterization of the intersection of $\mathcal{D}_d^{(0)}$ and any closed convex hull of finitely many interior points of $\mathcal{D}_d$ which is completely contained in the interior of $\mathcal{D}_d$. One of the algorithms is used to determine the structure of $\mathcal{D}_2^{(0)}$ in a region considerably larger than previously possible (see [18] for previous results), and to settle two important questions on its topology: It is shown that $\mathcal{D}_2^{(0)}$ is disconnected and that the largest connected component has a nontrivial fundamental group (Section 4). The other algorithm is of particular interest from a theoretic point of view as it is the first algorithm characterizing $\mathcal{D}_d^{(0)}$ in a given convex polyhedron which terminates for all inputs.

SRS are closely related to two important notions of numeration systems. Indeed, as pointed out in [1,8], SRS form a generalization of $\beta$-expansions (see [6,15,17]) and canonical number systems (CNS) (see [9,13,16] and [12, Section 4.1]).

For a non-integral real number $\beta > 1$ every $\gamma \in [0, \infty)$ can be represented uniquely in the form $\gamma = a_m\beta^m + a_{m-1}\beta^{m-1} + \ldots$ with $m \in \mathbb{Z}$ and $a_i \in \mathcal{A} := \{0, \ldots, \lfloor\beta\rfloor\}$ (*set of digits*) for all $m \geq i \in \mathbb{Z}$, such that $0 \leq \gamma - \sum_{i=k}^{m} a_i\beta^i < \beta^k$ holds for all $m \geq k \in \mathbb{Z}$ (*greedy expansion* of $\gamma$ with respect to $\beta$). Let $\operatorname{Fin}(\beta)$ be the set of all $\gamma \in [0,1)$ having finite greedy expansion with respect to $\beta$. Then $\operatorname{Fin}(\beta) \subseteq \mathbb{Z}[\frac{1}{\beta}] \cap [0,1)$. If the other inclusion is also true, $\beta$ is said to have *property (F)*. It is known that property (F) can only hold if $\beta$ is an algebraic integer and therefore has a minimal polynomial $X^d + a_{d-1}X^{d-1} + \cdots + a_1X + a_0 \in \mathbb{Z}[X]$ which can be written as $(X-\beta)(X^{d-1} + r_{d-2}X^{d-2} + \ldots + r_1X + r_0)$. On the other hand, not all algebraic integers have property (F). Indeed, $\beta$ has property (F) iff $(r_0, \ldots, r_{d-2}) \in \mathcal{D}_{d-1}^{(0)}$ [8].

[a] Note that this definition of SRS slightly differs from the original one in [1]. The SRS there are exactly those SRS (in the notion of the present paper) which have the finiteness property.
[b] From now on a real vector $\mathbf{r}$ and its associated SRS $\tau_{\mathbf{r}}$ shall be identified in terms of properties.

A similar relation can be shown for CNS. Let $P(X) = X^d + p_{d-1}X^{d-1} + \ldots + p_1X + p_0 \in \mathbb{Z}[X]$, $\mathcal{R} := \mathbb{Z}[X]/P(X)\mathbb{Z}[X]$, $\mathcal{N} := \{0, \ldots, |p_0| - 1\}$ and $x := X + P(X)\mathbb{Z}[X] \in \mathcal{R}$. Then $(P, \mathcal{N})$ is called a *canonical number system*, $P$ a *CNS polynomial* and $\mathcal{N}$ the *set of digits* iff every non-zero element $A(x) \in \mathcal{R}$ can be represented uniquely in the form $A(x) = a_mx^m + a_{m-1}x^{m-1} + \ldots + a_1x + a_0$ with $m \in \mathbb{N}_0$, $a_i \in \mathcal{N}$ for all $i \in \{0, \ldots, m\}$ and $a_m \neq 0$. Then $P$ is a CNS polynomial iff $(\frac{1}{p_0}, \frac{p_{d-1}}{p_0}, \ldots, \frac{p_2}{p_0}, \frac{p_1}{p_0}) \in \mathcal{D}_d^{(0)}$ [1].

We proceed with some technical preliminaries. For $n \in \mathbb{N}$, $\pi = (\mathbf{a}_1, \ldots \mathbf{a}_n) \in (\mathbb{Z}^d)^n$ is called a *cycle of* $\mathbf{r}$ (or $\tau_{\mathbf{r}}$, see footnote b) iff for all $i \in \{1, \ldots, n\}$ it holds that $\tau_{\mathbf{r}}(\mathbf{a}_i) = \mathbf{a}_{i \bmod n+1}$, a *cycle* iff there is a vector $\mathbf{r} \in \mathbb{R}^d$ for which $\pi$ is a cycle of $\mathbf{r}$, and *nontrivial* iff $\pi \neq (\mathbf{0})$, the *trivial cycle.* The *set of all cycles in* $\mathbb{Z}^d$ shall be denoted by $\mathcal{C}_d$. Let $P(\pi) := \{\mathbf{r} \in \mathbb{R}^d \mid \pi \text{ cycle of } \mathbf{r}\}$, the *associated polyhedron of* $\pi$ or - if $\pi$ is a nontrivial cycle - the *cutout polyhedron of* $\pi$. It follows from Lemma 2.2 below that $P(\pi)$ is either empty or the intersection of finitely many half-open “strips” and therefore it does in fact always form a - possibly degenerate - convex polyhedron. It is clear that

$$\mathcal{D}_d^{(0)} = \mathcal{D}_d \setminus \bigcup_{\pi \neq (\mathbf{0})} P(\pi) \tag{1.4}$$

which provides a method to “cut out” regions (the cutout polyhedra) from $\mathcal{D}_d$ [1].

The first algorithm described in this paper relies on a similar idea but calculates polyhedra which either belong entirely to $\mathcal{D}_d^{(0)}$ or have empty intersection with it. These polyhedra do not come from cycles but from so-called sets of witnesses which are used in what is known as Brunotte’s algorithm ([1, Theorem 5.1]). A set $V \subseteq \mathbb{Z}^d$ is called a *set of witnesses for* $\mathbf{r}$ iff it is stable under $\tau_{\mathbf{r}}$ and $\tau_{\mathbf{r}}^\star := -\tau_{\mathbf{r}} \circ (-\operatorname{id}_{\mathbb{Z}^d})$ and contains a generating set of the group $(\mathbb{Z}^d, +)$ which is closed under taking inverses. Every such set of witnesses has the decisive property

$$\mathbf{r} \in \mathcal{D}_d^{(0)} \Leftrightarrow \forall\, \mathbf{a} \in V : \exists\, n \in \mathbb{N} : \tau_{\mathbf{r}}^n(\mathbf{a}) = \mathbf{0}. \tag{1.5}$$

In the case of a finite set of witnesses this provides a method to decide whether or not a given parameter $\mathbf{r}$ belongs to $\mathcal{D}_d^{(0)}$. This is what Brunotte’s algorithm does for any given parameter $\mathbf{r}$ in the interior of $\mathcal{D}_d$ - it calculates a finite set of witnesses. It shall be denoted by $V_{\mathbf{r}}$ - the *set of witnesses associated with* $\mathbf{r}$ - and can be calculated using the following iteration ($V_{\mathbf{r}}$ can and shall be defined not only if $\mathbf{r}$ is in the interior of $\mathcal{D}_d$ but for any $\mathbf{r} \in \mathbb{R}^d$):

$$\begin{aligned} V_0 &:= \{(\pm 1, 0, \ldots, 0), \ldots, (0, \ldots, 0, \pm 1)\} \\ \forall\, n \in \mathbb{N} : \; V_n &:= V_{n-1} \cup \tau_{\mathbf{r}}(V_{n-1}) \cup \tau_{\mathbf{r}}^\star(V_{n-1}) \\ V_{\mathbf{r}} &:= \bigcup_{n \in \mathbb{N}_0} V_n \end{aligned} \tag{1.6}$$

If $\mathbf{r}$ is an element of the interior of $\mathcal{D}_d$ the mappings $\tau_{\mathbf{r}}$ and $\tau_{\mathbf{r}}^\star$ are contractive on some $\mathbb{Z}^d \setminus B_R(\mathbf{0})$ (where $B_R(\mathbf{0})$ denotes the $R$-ball around $\mathbf{0}$ for some $R \in \mathbb{R}$). Therefore the above iteration becomes stationary eventually [1]. Let $\Pi_{\mathbf{r}}$ - the *graph of witnesses associated with* $\mathbf{r}$ - denote the edge-colored multidigraph with vertex set $V_{\mathbf{r}}$ and an edge of color 1 from a vertex $\mathbf{a}$ to a vertex $\mathbf{b}$ iff $\tau_{\mathbf{r}}(\mathbf{a}) = \mathbf{b}$ and an edge of color 2 from $\mathbf{a}$ to $\mathbf{b}$ iff $\tau_{\mathbf{r}}^\star(\mathbf{a}) = \mathbf{b}$. If $E_1$ is the set of all edges (ordered pairs) of color 1 and $E_2$ the set of all edges of color 2, then the graph $\Pi_{\mathbf{r}}$ is completely characterized by the pair $(E_1, E_2) \in \mathcal{P}((\mathbb{Z}^d)^2)^2$ (as there are no isolated vertices) and thus the graph and the pair can be identified. For any such graph $\Pi = (E_1, E_2) \in \mathcal{P}((\mathbb{Z}^d)^2)^2$ let - just as for cycles - $P(\Pi) := \left\{\mathbf{r} \in \mathbb{R}^d \mid \forall\, (\mathbf{a}, \mathbf{b}) \in E_1 : \tau_{\mathbf{r}}(\mathbf{a}) = \mathbf{b} \,\wedge\, \forall\, (\mathbf{a}, \mathbf{b}) \in E_2 : \tau_{\mathbf{r}}^\star(\mathbf{a}) = \mathbf{b}\right\}$ and $P_{\mathbf{r}} := P(\Pi_{\mathbf{r}})$. As before one can use Lemma 2.2 from below to see that $P_{\mathbf{r}}$ is a convex polyhedron if $V_{\mathbf{r}}$ is finite, which is the case if $\mathbf{r} \in \operatorname{int}(\mathcal{D}_d)$. Furthermore $\mathcal{D}_d^{(0)}$ is the disjoint union of those $P_{\mathbf{r}}$ the corresponding parameters $\mathbf{r}$ of which belong to $\mathcal{D}_d^{(0)}$ (Proposition 3.1), which is essentially the heart of the first algorithm.

The second algorithm relies on a certain refinement of the decomposition of $\mathcal{D}_d^{(0)}$ provided by the set of the $P_{\mathbf{r}}$s, which can be computed much faster. However, this speedup comes with the same price Brunotte's algorithm for regions has to pay: the uncertainty whether or not the algorithm will terminate. But as with Brunotte's algorithm this limitation is only of theoretic interest and presents no difficulties in practice.

In addition to the regions of $\mathcal{D}_2^{(0)}$ which are settled algorithmically in Section 4, several infinite families of cutout polygons are deduced in Section 5. These families form a chain leading from $(1, -1)$ to $(1, 2)$ and tend to either of the two "critical points" (cf. [1]) $(1, 0)$ and $(1, 1)$.

## 2. Preliminary Results

**Definition 2.1.** *A* strip $S \subseteq \mathbb{R}^d$ *is the intersection of two parallel oppositely oriented half-spaces, or* $\mathbb{R}^d$ *itself. The empty set and the whole space* $\mathbb{R}^d$ *are considered* degenerate *and all others* nondegenerate *strips. For nondegenerate strips the attributes* open, half-open, *and* closed *shall indicate belonging of the hyperplanes bounding the strip. The* width *of a strip is the normal distance of these hyperplanes if it is nondegenerate or* $-\infty$ *or* $\infty$ *if the strip is the empty set or the whole space.*

*A set* $P \subseteq \mathbb{R}^d$ *is called* (convex) polyhedron *iff it is the intersection of finitely many half-spaces or* $\mathbb{R}^d$ *itself. A polyhedron is considered* nondegenerate *iff it has positive and finite Lebesgue measure and* degenerate *otherwise. The set of all polyhedra in* $\mathbb{R}^d$ *shall be denoted by* $\mathcal{P}_d$.

Every strip having positive width can be represented in one of the four ways

$$\left\{\mathbf{r} \in \mathbb{R}^d \mid 0 \left\{\begin{matrix} < \\ \leq \end{matrix}\right\} \mathbf{a}\mathbf{r} + b \left\{\begin{matrix} < \\ \leq \end{matrix}\right\} 1\right\} \qquad (\mathbf{a} \in \mathbb{R}^d,\, b \in \mathbb{R}), \tag{2.1}$$

where $\mathbf{a}$ is normal to the strip's bounding hyperplanes and $\frac{1}{\|\mathbf{a}\|}$ is the strip's width.

**Lemma 2.2.** *Let* $\mathbf{a} = (a_1, \ldots, a_d) \in \mathbb{Z}^d$, $\mathbf{b} = (b_1, \ldots, b_d) \in \mathbb{Z}^d$, *and* $\mathbf{r} \in \mathbb{R}^d$. *Then*

*(i)* $\left\{\mathbf{r} \in \mathbb{R}^d \mid \tau_{\mathbf{r}}(\mathbf{a}) = \mathbf{b}\right\} =$
$\left\{\mathbf{r} \in \mathbb{R}^d \mid (a_2, \ldots, a_d) = (b_1, \ldots, b_{d-1}) \wedge 0 \le \mathbf{ra} + b_d < 1\right\}$

*(ii)* $\left\{\mathbf{r} \in \mathbb{R}^d \mid \tau^\star_{\mathbf{r}}(\mathbf{a}) = \mathbf{b}\right\} =$
$\left\{\mathbf{r} \in \mathbb{R}^d \mid (a_2, \ldots, a_d) = (b_1, \ldots, b_{d-1}) \wedge 0 \le -\mathbf{ra} - b_d < 1\right\}$

*(iii)* $\left\{\mathbf{s} \in \mathbb{R}^d \mid \tau_{\mathbf{s}}(\mathbf{a}) = \tau_{\mathbf{r}}(\mathbf{a}) \wedge \tau^\star_{\mathbf{s}}(\mathbf{a}) = \tau^\star_{\mathbf{r}}(\mathbf{a})\right\} =$

$$\begin{cases} \left\{\mathbf{s} \in \mathbb{R}^d \mid \mathbf{sa} - \mathbf{ra} = 0\right\} & \text{if } \mathbf{ra} \in \mathbb{Z} \\ \left\{\mathbf{s} \in \mathbb{R}^d \mid 0 < \mathbf{sa} - \lfloor \mathbf{ra} \rfloor < 1\right\} & \text{if } \mathbf{ra} \notin \mathbb{Z} \end{cases} \quad \textit{(hyperplane or open strip)}$$

**Proof.** (i) and (ii) follow directly from the definitions of $\tau_{\mathbf{r}}$ and $\tau^\star_{\mathbf{r}}$.

For the proof of (iii) let $M := \left\{\mathbf{s} \in \mathbb{R}^d \mid \tau_{\mathbf{s}}(\mathbf{a}) = \tau_{\mathbf{r}}(\mathbf{a}) \wedge \tau^\star_{\mathbf{s}}(\mathbf{a}) = \tau^\star_{\mathbf{r}}(\mathbf{a})\right\}$. Then (i) and (ii) imply that $M = \left\{\mathbf{s} \in \mathbb{R}^d \mid 0 \le \mathbf{sa} - \lfloor \mathbf{ra} \rfloor < 1 \wedge 0 \le -\mathbf{sa} - \lfloor -\mathbf{ra} \rfloor < 1\right\}$. If $\mathbf{ra} \in \mathbb{Z}$ then $-\lfloor -\mathbf{ra} \rfloor = \mathbf{ra}$ and therefore $M = \left\{\mathbf{s} \in \mathbb{R}^d \mid \mathbf{sa} - \mathbf{ra} = 0\right\}$. If $\mathbf{ra} \notin \mathbb{Z}$ then $-\lfloor -\mathbf{ra} \rfloor = \lfloor \mathbf{ra} \rfloor + 1$ and therefore $M = \left\{\mathbf{s} \in \mathbb{R}^d \mid 0 < \mathbf{sa} - \lfloor \mathbf{ra} \rfloor < 1\right\}$. □

**Lemma 2.3.** *Let* $\mathbf{r} \in \operatorname{int}(\mathcal{D}_d)$. *Then* $P_{\mathbf{r}}$ *is the intersection of a nondegenerate, open polyhedron and an affine subspace of* $\mathbb{R}^d$.

**Proof.** Lemma 2.2 implies that $P_{\mathbf{r}}$ is the intersection of finitely many hyperplanes and finitely many open strips and thus is the intersection of an open polyhedron and an affine subspace of $\mathbb{R}^d$. Furthermore it is non-empty as $\mathbf{r} \in P_{\mathbf{r}}$ and it is bounded as $\{(\pm 1, 0, \ldots, 0), \ldots, (0, \ldots, 0, \pm 1)\} \subseteq V_{\mathbf{r}}$ and therefore $P_{\mathbf{r}} \subseteq \lfloor \mathbf{r} \rfloor + [0,1]^d$. □

## 3. Two Characterization Algorithms for Regions of $\mathcal{D}_d^{(0)}$

Throughout the following section let $k \in \mathbb{N}$, $(\mathbf{r}_1, \ldots, \mathbf{r}_k) \in \operatorname{int}(\mathcal{D}_d)^k$ and $H := \operatorname{conv}(\{\mathbf{r}_1, \ldots, \mathbf{r}_k\})$ such that $H \subset \operatorname{int}(\mathcal{D}_d)$.

**Algorithm 1** Determination of $\operatorname{conv}(\{\mathbf{r}_1, \ldots, \mathbf{r}_k\}) \cap \mathcal{D}_d^{(0)}$

**Input:** $(\mathbf{r}_1, \ldots, \mathbf{r}_k) \in \operatorname{int}(\mathcal{D}_d)^k$ such that $\operatorname{conv}(\{\mathbf{r}_1, \ldots, \mathbf{r}_k\}) \subset \operatorname{int}(\mathcal{D}_d)$.
**Output:** $\mathcal{P} \subseteq \mathcal{P}_d$ with $\operatorname{conv}(\{\mathbf{r}_1, \ldots, \mathbf{r}_k\}) \cap \mathcal{D}_d^{(0)} = \bigcup \mathcal{P}$ disjoint.

$H \leftarrow \operatorname{conv}(\{\mathbf{r}_1, \ldots, \mathbf{r}_k\})$
$\mathcal{P} \leftarrow \emptyset$
**while** $H \setminus \bigcup \mathcal{P} \neq \emptyset$ **do**
  select $\mathbf{r} \in H \setminus \bigcup \mathcal{P}$
  $\mathcal{P} \leftarrow \mathcal{P} \cup \{H \cap P_{\mathbf{r}}\}$
  **if** $\mathbf{r} \in \mathcal{D}_d^{(0)}$ **then** {use Brunotte's algorithm}
    $\operatorname{fin}_{H \cap P_{\mathbf{r}}} \leftarrow$ **true**
  **else**
    $\operatorname{fin}_{H \cap P_{\mathbf{r}}} \leftarrow$ **false**
  **end if**
**end while**
**return** $\{P \in \mathcal{P} \mid \operatorname{fin}_P = \textbf{true}\}$

Algorithm 1 is a straightforward application of the following

**Proposition 3.1.** *Let* $\mathbf{r} \in \mathbb{R}^d$. *Then*

*(i)* $\mathbf{r} \in P_{\mathbf{r}}$, *and* $\mathbf{r} \in \mathcal{D}_d^{(0)} \Leftrightarrow P_{\mathbf{r}} \subseteq \mathcal{D}_d^{(0)}$
*(ii)* $\mathcal{D}_d^{(0)} = \bigcup \left\{ P_{\mathbf{r}} \mid \mathbf{r} \in \mathcal{D}_d^{(0)} \right\}$ *and this union is disjoint*

**Proof.** Follows directly from the definition of $P_{\mathbf{r}}$. □

Of course the question arises whether the while loop actually terminates, which is equivalent to the possibility of exhausting $H$ by finitely many $P_{\mathbf{r}}$s.

**Theorem 3.2.** *Algorithm 1 terminates for all inputs.*

**Proof.** In [1, Section 4] it is shown that $\tau_{\mathbf{r}}$ can be represented by the transposed companion matrix $R_{\mathbf{r}}$ of the polynomial $\chi_{\mathbf{r}}(X) = X^d + r_d X^{d-1} + \cdots + r_2 X + r_1$, where $\mathbf{r} = (r_1, \ldots, r_d) \in \mathbb{R}^d$. For every $\mathbf{a} \in \mathbb{Z}^d$ there is a $\mathbf{v_r} = (0, \ldots, 0, v_{\mathbf{r}}) \in \mathbb{R}^d$ with $0 \leq v_{\mathbf{r}} < 1$ such that $\tau_{\mathbf{r}}(\mathbf{a}) = R_{\mathbf{r}}\mathbf{a} + \mathbf{v_r}$.

Furthermore formula (3.2) of [14] is used to show that for every $\mathbf{r} \in \mathbb{R}^d$ and $\tilde{\rho} > \rho(R_{\mathbf{r}})$ (the spectral radius of $R_{\mathbf{r}}$) there is a norm $\|\cdot\|_{\mathbf{r},\tilde{\rho}}$ such that $\|R_{\mathbf{r}}\mathbf{a}\|_{\mathbf{r},\tilde{\rho}} \leq \tilde{\rho}\,\|\mathbf{a}\|_{\mathbf{r},\tilde{\rho}}$ for all $\mathbf{a} \in \mathbb{R}^d$ and $\|R_{\mathbf{r}}\mathbf{a}\|_{\mathbf{r},\tilde{\rho}} = \tilde{\rho}\,\|\mathbf{a}\|_{\mathbf{r},\tilde{\rho}}$ iff $\mathbf{a} = \mathbf{0}$. As the function which maps $\mathbf{s} \in \mathbb{R}^d$ to $\max\left\{ \frac{\|R_{\mathbf{s}}\mathbf{a}\|_{\mathbf{r},\tilde{\rho}}}{\|\mathbf{a}\|_{\mathbf{r},\tilde{\rho}}} \mid \mathbf{a} \in \mathbb{R}^d \wedge \|\mathbf{a}\|_{\mathbf{r},\tilde{\rho}} = 1 \right\}$ is continuous, there is an open neighborhood $B_{\mathbf{r},\tilde{\rho}}$ of $\mathbf{r}$ such that $\|R_{\mathbf{s}}\mathbf{a}\|_{\mathbf{r},\tilde{\rho}} \leq \tilde{\rho}\,\|\mathbf{a}\|_{\mathbf{r},\tilde{\rho}}$ for every $\mathbf{s} \in B_{\mathbf{r},\tilde{\rho}}$ and $\mathbf{a} \in \mathbb{R}^d$.

Since $\mathbf{r}$ is an interior point of $\mathcal{D}_d$ iff $\rho(R_{\mathbf{r}}) < 1$ ([1, Lemma 4.1 and Lemma 4.3]) we get that for every $\mathbf{r} \in H \subset \operatorname{int}(\mathcal{D}_d)$ there is a $\rho(R_{\mathbf{r}}) < \tilde{\rho_{\mathbf{r}}} < 1$ such that $\|\tau_{\mathbf{s}}(\mathbf{a})\|_{\mathbf{r},\tilde{\rho_{\mathbf{r}}}} < \|\mathbf{a}\|_{\mathbf{r},\tilde{\rho_{\mathbf{r}}}}$ for every $\mathbf{s} \in B_{\mathbf{r},\tilde{\rho_{\mathbf{r}}}}$ and $\mathbf{a} \in \mathbb{Z}^d$ with $\|\mathbf{a}\|_{\mathbf{r},\tilde{\rho_{\mathbf{r}}}} > \rho_{\mathbf{r}} := \frac{\max\left\{\|(0,\ldots,0,v)\|_{\mathbf{r},\tilde{\rho_{\mathbf{r}}}} | v \in [0,1]\right\}}{1-\tilde{\rho_{\mathbf{r}}}}$. This implies that $V_{\mathbf{s}} \subseteq \left\{ \mathbf{a} \in \mathbb{Z}^d \mid \|\mathbf{a}\|_{\mathbf{r},\tilde{\rho_{\mathbf{r}}}} \leq \rho_{\mathbf{r}} \right\}$ for every $\mathbf{s} \in B_{\mathbf{r},\tilde{\rho_{\mathbf{r}}}}$. The set $\{B_{\mathbf{r},\tilde{\rho_{\mathbf{r}}}} \mid \mathbf{r} \in H\}$ is an open cover of $H$ and since $H$ is compact there exists a finite subcover. Thus $\bigcup_{\mathbf{r} \in H} V_{\mathbf{r}}$ is finite and since only finitely many graphs can be defined on a finite number of vertices the sets $\{\Pi_{\mathbf{r}} \mid \mathbf{r} \in H\}$ and $\{P_{\mathbf{r}} \mid \mathbf{r} \in H\}$ are also finite. □

Algorithm 1 computes a decomposition of $H$ into finitely many disjoint polyhedra (from which it selects those which are contained in $\mathcal{D}_d^{(0)}$ in the final step). Algorithm 2 - which will be treated in the following - uses any finite superset $V$ (which has to be fixed initially) of $\tilde{V}_H := \bigcup_{\mathbf{r} \in H} V_{\mathbf{r}}$ to compute a refinement of this decomposition (cf. Fig. 1). $\tilde{V}_H$ itself is a finite set according to Theorem 3.2 and can be computed by Algorithm 1. Though $\tilde{V}_H$ would be the optimal choice for $V$ its determination by Algorithm 1 would of course be pointless. But at least for some choices of $H$ another finite superset of $\tilde{V}_H$ can be calculated efficiently using Brunotte's algorithm for regions [1]. It calculates a common set of witnesses for all $\mathbf{r} \in H$ essentially using the same iteration as for single parameters (cf. formula (1.6)) but with the slightly modified functions $\overline{\tau}_H : \mathcal{P}(\mathbb{Z}^d) \to \mathcal{P}(\mathbb{Z}^d)$, where

$$\overline{\tau}_H(V) = \{\tau_{\mathbf{r}}(\mathbf{a}) \mid \mathbf{r} \in H \wedge \mathbf{a} \in V\}\,, \tag{3.1}$$

and $\overline{\tau}_H^\star := -\overline{\tau}_H \circ (-\operatorname{id}_{\mathcal{P}(\mathbb{Z}^d)})$. A possible choice for $V_H$ can then be found using the following iteration:

$$\begin{aligned} V_0 &:= \{(\pm 1, 0, \ldots, 0), \ldots, (0, \ldots, 0, \pm 1)\} \\ \forall\, n \in \mathbb{N}:\ V_n &:= V_{n-1} \cup \overline{\tau}_H(V_{n-1}) \cup \overline{\tau}_H^\star(V_{n-1}) \\ V_H &:= \bigcup_{n\in\mathbb{N}_0} V_n \end{aligned} \tag{3.2}$$

Unfortunately the set $V_H$ found in this way need not always be finite even if $H$ is contained in the interior of $\mathcal{D}_d$. However in practice this causes no troubles as one just has to start with a smaller $H$ if the iteration does not become stationary, which can be decided heuristically by comparing the size (i.e. number of elements or maximum of its elements' norms) of $V_n$ to the sizes of the $V_{\mathbf{r}}$s for all vertices $\mathbf{r}$ of the polyhedron $H$. If $V_n$ is considerably larger than all of the $V_{\mathbf{r}}$s the iteration probably will not terminate.

From now on let $V \subseteq \mathbb{Z}^d$ fix any finite superset of $\tilde{V}_H$ and consider the following equivalence relation:

$$\sim := \left\{(\mathbf{r}_1, \mathbf{r}_2) \in \mathbb{R}^2 \mid \forall\, \mathbf{a} \in V : \tau_{\mathbf{r}_1}(\mathbf{a}) = \tau_{\mathbf{r}_2}(\mathbf{a}) \wedge \tau^\star_{\mathbf{r}_1}(\mathbf{a}) = \tau^\star_{\mathbf{r}_2}(\mathbf{a})\right\} \tag{3.3}$$

Then the set $H/_\sim := \left\{[\mathbf{r}]_\sim \cap H \mid \mathbf{r} \in \mathbb{R}^d\right\}$ is a refinement of the decomposition of $H$ calculated by Algorithm 1. If $\mathcal{R} \subseteq H$ is any system of representatives of $H/_\sim$ then the intersection of $\mathcal{D}_d^{(0)}$ and $H$ is given by the finite disjoint union $\bigcup\left\{[\mathbf{r}]_\sim \cap H \mid \mathbf{r} \in \mathcal{R} \cap \mathcal{D}_d^{(0)}\right\}$. In order to determine the complete list of equivalence classes (Proposition 3.5 below) the notion of *face lattices* of (convex) polyhedra proves useful (an adapted version of [7, Chapter 3] is used to cover degenerate polyhedra which will be needed in Section 5).

**Definition 3.3.** *A* face *of a polyhedron $P \in \mathcal{P}_d$ is any intersection of $\overline{P}$ with a closed half-space such that the interior of $P$ (with respect to the smallest affine subspace of $\mathbb{R}^d$ containing $P$) and the boundary of the half-space are disjoint. In addition $\emptyset$ and $\mathbb{R}^d$ shall be considered faces if $P = \mathbb{R}^d$. The set of faces of $P$ shall be denoted by $\mathcal{F}(P)$.*

*The* face lattice *of $P$ is the set of faces $\mathcal{F}(P)$ of $P$ together with the partial order given by set inclusion.*

*For a face $F$ of $P$ let $F^\circ$ denote the set difference of $F$ and the union of all faces of $P$ that are less than $F$ (in the face lattice of $P$). Any $F^\circ$ where $F \in \mathcal{F}(P)$ shall be referred to as* open face *of $P$ and consequently $\mathcal{F}^\circ(P) := \{F^\circ \mid F \in \mathcal{F}(P)\}$ as the* set of open faces *of $P$.*

Note that any closed polyhedron is the disjoint union of its open faces. Also note that open faces need not to be open sets (e.g. the singletons containing the vertices of a polyhedron are among its open faces). Furthermore we need the following technical

**Definition 3.4.** *For a tuple* $(\mathbf{a}, b) \in \mathbb{R}^d \times \mathbb{R}$ *let* $\mathrm{P}(\mathbf{a}, b) := \left\{\mathbf{r} \in \mathbb{R}^d \mid \mathbf{a}\mathbf{r} + b = 0\right\}$. *A hyperplane* $P \subseteq \mathbb{R}^d$ *is called* integer *if there is a tuple* $(\mathbf{a}, b) \in \mathbb{Z}^d \times \mathbb{Z}$ *such that* $P = \mathrm{P}(\mathbf{a}, b)$. *Any such tuple shall then be denoted as* generator of $P$. *The unique generator which satisfies that the first nonzero entry of* $\mathbf{a}$ *is positive and that the greatest common divisor of the entries of* $\mathbf{a}$ *and* $b$ *is* $1$ *is the* canonical generator of $P$ *and shall be denoted by* $\mathrm{CG}(P) = (\mathrm{CG}_1(P), \mathrm{CG}_2(P))$. *The first entry* $\mathrm{CG}_1(P)$ *of the canonical generator is the* canonical normal vector of $P$.

**Proposition 3.5.** *For all* $\mathbf{a} \in V$ *let* $B_{\mathbf{a}} := \{-\mathbf{a}\mathbf{r}_i \mid i \in \{1, \ldots, k\}\}$ *and* $\mathcal{G} := \{\mathrm{CG}(\mathrm{P}(\mathbf{a}, b)) \mid \mathbf{a} \in V \setminus \{\mathbf{0}\} \wedge b \in \{\lfloor \min(B_{\mathbf{a}}) \rfloor, \ldots, \lceil \max(B_{\mathbf{a}}) \rceil\}\}$. *Furthermore let* $\phi : \mathbb{R}^d \to \{-1, 0, 1\}^{|\mathcal{G}|}$ *where* $\phi(\mathbf{r}) = (\mathrm{sgn}(\mathbf{a}\mathbf{r} + b))_{(\mathbf{a},b)\in\mathcal{G}}$ *and let* $\mathcal{P}$ *denote the set of all minimal nondegenerate polyhedra having non-empty intersection with* $H$ *which are the intersection of some selection of half-spaces from the set* $\left\{\left\{\mathbf{r} \in \mathbb{R}^d \mid \mathbf{a}\mathbf{r} + b \geq 0\right\} \mid (\mathbf{a}, b) \in \mathcal{G}\right\} \cup \left\{\left\{\mathbf{r} \in \mathbb{R}^d \mid -\mathbf{a}\mathbf{r} - b \geq 0\right\} \mid (\mathbf{a}, b) \in \mathcal{G}\right\}$. *Then* $H/_{\sim} = \left\{\phi^{-1}(S) \cap H \mid S \in \{-1, 0, 1\}^{|\mathcal{G}|}\right\} \setminus \{\emptyset\} = \left\{F \cap H \mid F \in \bigcup_{P\in\mathcal{P}} \mathcal{F}^{\circ}(P)\right\} \setminus \{\emptyset\}$.

**Proof.** Lemma 2.2 implies that $\Phi_1 : H/_{\sim} \to \{-1, 0, 1\}^{|\mathcal{G}|}$, where $\Phi_1([\mathbf{r}]_{\sim}) = \phi(\mathbf{r})$ is well-defined and injective. On the other hand it follows from the definitions of open faces and $\mathcal{P}$ that $\Phi_2 : \left\{F \cap H \mid F \in \bigcup_{P\in\mathcal{P}} \mathcal{F}^{\circ}(P)\right\} \setminus \{\emptyset\} \to \{-1, 0, 1\}^{|\mathcal{G}|}$ given by $\Phi_2(F) = (\mathrm{sgn}(\mathbf{a}\mathbf{v} + b))_{(\mathbf{a},b)\in\mathcal{G}}$ with $\mathbf{v} \in F$ is also well-defined and injective and $\Phi_1^{-1}(S) = \Phi_2^{-1}(S)$ for all $S \in \{-1, 0, 1\}^{|\mathcal{G}|}$ which proves the statement. □

**Remark 3.6.** Proposition 3.5 gives a geometric interpretation of the equivalence classes of $\sim$. The hyperplanes which are generated by the elements of $\mathcal{G}$ cut $\mathbb{R}^d$ into pieces of polyhedral shape and the set of all (nonempty) open faces of these polyhedra is exactly the set of equivalence classes of $\sim$. The use of canonical generators eliminates redundant hyperplanes, which is not needed in the proof but will speed up the process of actually finding the set of all open faces. If $d = 2$ this is not too difficult but one would probably approach the problem in reverse order than what Proposition 3.5 suggests. Instead of calculating the set $\mathcal{P}$ of polygons directly and then the set of open faces (singletons (vertices), open line segments (edges), and nondegenerate open polygons (interiors)) one can first find all vertices by pairwise intersection of the given lines, then the edges (pair of distinct vertices that lie on a common line with no other vertex lying in between), and at last the interiors (use any algorithm to find the graph theoretic faces of a planar embedding of a graph).

In higher dimensions one could use the cylindrical algebraic decomposition algorithm [5] which, for a given set of polynomials in $\mathbb{R}[X_1, \ldots, X_d]$, finds a decomposition of $\mathbb{R}^d$ into regions on which each polynomial has constant sign.

**Remark 3.7.** The set $H/_{\sim}$ of equivalence classes is also useful when calculating $\overline{\tau}_H(V)$ for some finite $V \subseteq \mathbb{Z}^d$ (see formula (3.2)). It follows from the definition of $\overline{\tau}_H$ that $\overline{\tau}_H(V) = \bigcup_{\mathbf{a}\in V} \overline{\tau}_H(\{\mathbf{a}\})$ and for any $\mathbf{a} \in \mathbb{Z}^d$ one gets that $\overline{\tau}_H(\{\mathbf{a}\}) = \{\tau_{\mathbf{r}}(\mathbf{a}) \mid [\mathbf{r}]_{\sim} \in H/_{\sim}\}$ where $\sim := \left\{(\mathbf{r}_1, \mathbf{r}_2) \in \mathbb{R}^2 \mid \tau_{\mathbf{r}_1}(\mathbf{a}) = \tau_{\mathbf{r}_2}(\mathbf{a}) \wedge \tau^{\star}_{\mathbf{r}_1}(\mathbf{a}) = \tau^{\star}_{\mathbf{r}_2}(\mathbf{a})\right\}$.

If the set $H/_\sim$ of equivalence classes is known one could use Brunotte's algorithm to decide whether or not a given class belongs to $\mathcal{D}_d^{(0)}$. The definition of $\sim$ guarantees that the result will be the same for all parameters in the class. But instead of treating all classes independently two decisive optimizations can be made to speed up the process considerably.

If $[\mathbf{r}]_\sim \in H/_\sim$ is any class then $V_\mathbf{r} \subseteq V$ and in all situations of practical relevance (where $V$ has been found with Brunotte's algorithm for regions) $V$ probably will not be much larger than $V_\mathbf{r}$. So instead of calculating $V_\mathbf{r}$ and checking it for nontrivial cycles one can just check the similar superset $V$.

The second optimization relies on the fact that the graph defined by $\tau_\mathbf{r}$ on $V$ can only change at specific vertices if the parameter is changed to $\mathbf{s}$ where $[\mathbf{s}]_\sim \in H/_\sim$ is any class that is adjacent to $[\mathbf{r}]_\sim$ and two classes are considered *adjacent* if they are distinct and their topological boundaries intersect. The following proposition describes on which nodes the graph on $V$ needs to be updated (at most) in this situation. If both classes have a positive distance from the boundary of $H$, the set $M$ of these nodes consists of those elements of $V$ which are integer multiples of the canonical normal vectors of any hyperplane containing the intersection of the boundaries of $[\mathbf{r}]_\sim$ and $[\mathbf{s}]_\sim$ and any (d-1)-dimensional class in the intersection of the "closed neighborhoods" of $[\mathbf{r}]_\sim$ and $[\mathbf{s}]_\sim$.

**Proposition 3.8.** *Let* $[\mathbf{r}]_\sim \in H/_\sim$ *and* $[\mathbf{s}]_\sim \in H/_\sim$ *be adjacent and for any class* $C \in H/_\sim$ *let* $\overline{N(C)} := \{D \in H/_\sim \mid D \text{ adjacent to } C\} \cup \{C\}$ *(closed neighborhood of* $C$*). Furthermore let* $M := \{\mathbf{a} \in V \mid \exists\, b \in \mathbb{Z} : \partial[\mathbf{r}]_\sim \cap \partial[\mathbf{s}]_\sim \subseteq \mathrm{P}(\mathbf{a}, b)\}$*. Then*

*(i)* $\{\mathbf{a} \in V \mid \tau_\mathbf{r}(\mathbf{a}) \neq \tau_\mathbf{s}(\mathbf{a})\} \subseteq M$

*(ii)* $\partial[\mathbf{r}]_\sim \cap \partial H = \emptyset \;\wedge\; \partial[\mathbf{s}]_\sim \cap \partial H = \emptyset \Rightarrow$

$$\begin{aligned} M = V \cap \mathbb{Z}\{\mathrm{CG}_1(\mathrm{span}_\mathbb{R}([\mathbf{t}]_\sim - \mathbf{t}) + \mathbf{t}) \mid\ & [\mathbf{t}]_\sim \in \overline{N([\mathbf{r}]_\sim)} \cap \overline{N([\mathbf{s}]_\sim)} \wedge \\ & \dim_\mathbb{R}(\mathrm{span}_\mathbb{R}([\mathbf{t}]_\sim - \mathbf{t})) = d-1 \wedge \\ & \partial[\mathbf{r}]_\sim \cap \partial[\mathbf{s}]_\sim \subseteq \mathrm{span}_\mathbb{R}([\mathbf{t}]_\sim - \mathbf{t}) + \mathbf{t}\} \end{aligned}$$

**Proof.** We say that a hyperplane *separates* two classes from $H/_\sim$ iff either one class is contained in the hyperplane while the other has empty intersection with it or each class has empty intersection with exactly one of the two open half-spaces $\mathbb{R}^d$ is divided into by the hyperplane. For every $\mathbf{a} \in V$ with $\tau_\mathbf{r}(\mathbf{a}) \neq \tau_\mathbf{s}(\mathbf{a})$ there is a $b \in \mathbb{Z}$ such that $\mathrm{P}(\mathbf{a}, b)$ separates $[\mathbf{r}]_\sim$ and $[\mathbf{s}]_\sim$ according to Proposition 3.5. And every hyperplane separating the two distinct but "touching" classes $[\mathbf{r}]_\sim$ and $[\mathbf{s}]_\sim$ has to contain the intersection of their topological boundaries which shows (i).

For the proof of (ii) assume that $[\mathbf{r}]_\sim$ and $[\mathbf{s}]_\sim$ both have a positive distance from the boundary of $H$ and let $\mathbf{a} \in V$ and $b \in \mathbb{Z}$ such that $\partial[\mathbf{r}]_\sim \cap \partial[\mathbf{s}]_\sim \subseteq \mathrm{P}(\mathbf{a}, b)$. Then it follows from Proposition 3.5 that there is a class $[\mathbf{t}]_\sim \in \overline{N([\mathbf{r}]_\sim)} \cap \overline{N([\mathbf{s}]_\sim)}$ satisfying that $\mathrm{P}(\mathbf{a}, b) = \mathrm{span}_\mathbb{R}([\mathbf{t}]_\sim - \mathbf{t}) + \mathbf{t}$ which implies that $\dim_\mathbb{R}(\mathrm{span}_\mathbb{R}([\mathbf{t}]_\sim - \mathbf{t})) = \dim_\mathbb{R}(\mathrm{P}(\mathbf{a}, b) - \mathbf{t}) = d-1$. Furthermore it follows from the definition of the canonical generator that $\mathbf{a} \in \mathbb{Z}\,\mathrm{CG}_1(\mathrm{P}(\mathbf{a}, b))$ which proves (ii). □

**Algorithm 2** Determination of $\operatorname{conv}(\{\mathbf{r}_1,\dots,\mathbf{r}_k\}) \cap \mathcal{D}_d^{(0)}$

**Input:** $(\mathbf{r}_1,\dots,\mathbf{r}_k) \in \operatorname{int}(\mathcal{D}_d)^k$ such that $\operatorname{conv}(\{\mathbf{r}_1,\dots,\mathbf{r}_k\}) \subset \operatorname{int}(\mathcal{D}_d)$, $\tilde{V}_{\operatorname{conv}(\{\mathbf{r}_1,\dots,\mathbf{r}_k\})} \subseteq V \subseteq \mathbb{Z}^d$ finite.

**Output:** $\mathcal{C} \subseteq \mathcal{C}_d$ with $\operatorname{conv}(\{\mathbf{r}_1,\dots,\mathbf{r}_k\}) \cap \mathcal{D}_d^{(0)} = \operatorname{conv}(\{\mathbf{r}_1,\dots,\mathbf{r}_k\}) \setminus \bigcup_{\pi\in\mathcal{C}} P(\pi)$.

$H \leftarrow \operatorname{conv}(\{\mathbf{r}_1,\dots,\mathbf{r}_k\})$
$\mathcal{C} \leftarrow \emptyset$
$G = (V(G), E(G)) \leftarrow (V, \emptyset)$ {edgeless digraph with vertex set $V$}
calculate $H/_\sim$ according to Proposition 3.5 and Remark 3.6
**for all** $C \in H/_\sim$ **do**
    $N_C \leftarrow \{D \in H/_\sim \mid D$ adjacent to $C\}$
    $B_C \leftarrow$ **false**
**end for**
**for all** $[\mathbf{r}]_\sim \in H/_\sim$ with $B_{[\mathbf{r}]_\sim} =$ **false** and $\partial[\mathbf{r}]_\sim \cap \partial H \neq \emptyset$ **do**
    **if** $\mathbf{r} \in \mathcal{D}_d^{(0)}$ **then** {search for cycles of $\mathbf{r}$ on $V$}
        $B_{[\mathbf{r}]_\sim} \leftarrow$ **true**
    **else**
        select $\pi$ nontrivial cycle of $\mathbf{r}$ on $V$
        $\mathcal{C} \leftarrow \mathcal{C} \cup \{\pi\}$
        **for all** $[\mathbf{s}]_\sim \in H/_\sim$ with $B_{[\mathbf{s}]_\sim} =$ **false** and $\mathbf{s} \in P(\pi)$ **do**
            $B_{[\mathbf{s}]_\sim} \leftarrow$ **true**
        **end for**
    **end if**
**end for**
**while** $\exists\, C \in H/_\sim : B_C =$ **false do**
    select $[\mathbf{r}]_\sim \in H/_\sim$ with $B_{[\mathbf{r}]_\sim} =$ **false**
    $E(G) \leftarrow \{(\mathbf{a}, \tau_{\mathbf{r}}(\mathbf{a})) \mid \mathbf{a} \in V\}$
    $W \leftarrow V$
    **loop**
        **if** $\mathbf{r} \in \mathcal{D}_d^{(0)}$ **then** {search for cycles of $G$ starting at the vertices in $W$}
            $B_{[\mathbf{r}]_\sim} \leftarrow$ **true**
            **if** $\exists\, C \in N_{[\mathbf{r}]_\sim} : B_C =$ **false then**
                select $C \in N_{[\mathbf{r}]_\sim}$ with $B_C =$ **false**
                update $E(G)$ according to Proposition 3.8
                save the tails of the changed edges in $W$
                $[\mathbf{r}]_\sim \leftarrow C$
            **else**
                **break**
            **end if**
        **else**
            select $\pi$ nontrivial cycle of $G$
            $\mathcal{C} \leftarrow \mathcal{C} \cup \{\pi\}$
            **for all** $[\mathbf{s}]_\sim \in H/_\sim$ with $B_{[\mathbf{s}]_\sim} =$ **false** and $\mathbf{s} \in P(\pi)$ **do**
                $B_{[\mathbf{s}]_\sim} \leftarrow$ **true**
            **end for**
            **break**
        **end if**
    **end loop**
**end while**
**return** $\mathcal{C}$

Algorithm 2 computes a minimal set of cutout polyhedra (with respect to set inclusion but not necessarily cardinality) which describes $\mathcal{D}_d^{(0)}$ inside of $H$. After initialization of required variables (steps 1-8) the classes "touching" the boundary of $H$ are treated directly and one by one (steps 9-19). If a nontrivial cycle is found all classes contained in the corresponding cutout polyhedron are considered "done" (the associated $B$-flag is set to **true**) and the cycle is added to the output set $\mathcal{C}$.

After that the main part of the algorithm follows (steps 20-44). The classes are treated along walks in the graph defined on $H/_\sim$ by the adjacency relation. When a nontrivial cycle is found it is handled as before and a new walk begins, as it does when the walk reaches a dead end (i.e. if there are no neighbors yet to be treated). Any time a new walk starts all edges of $G$ have to be updated and checked for cycles which consumes much more time than updating and checking only those edges which are changed when going from one class to an adjacent one. In order to minimize the number of restarts it is crucial to make a good choice when selecting the next node (step 28). A Hamiltonian path would of course be an optimal but also costly choice. Instead the following heuristic turns out to be adequate: Of all possible neighbors of least dimension take the one (or one of those) which has the highest number of neighbors that are already treated. This way the walks tend to stay "compact" and will not cut the graph into too many pieces of pending vertices.

Fig. 1 illustrates the relation between the resulting decompositions of Algorithm 1 and Algorithm 2 and cutout polygons.

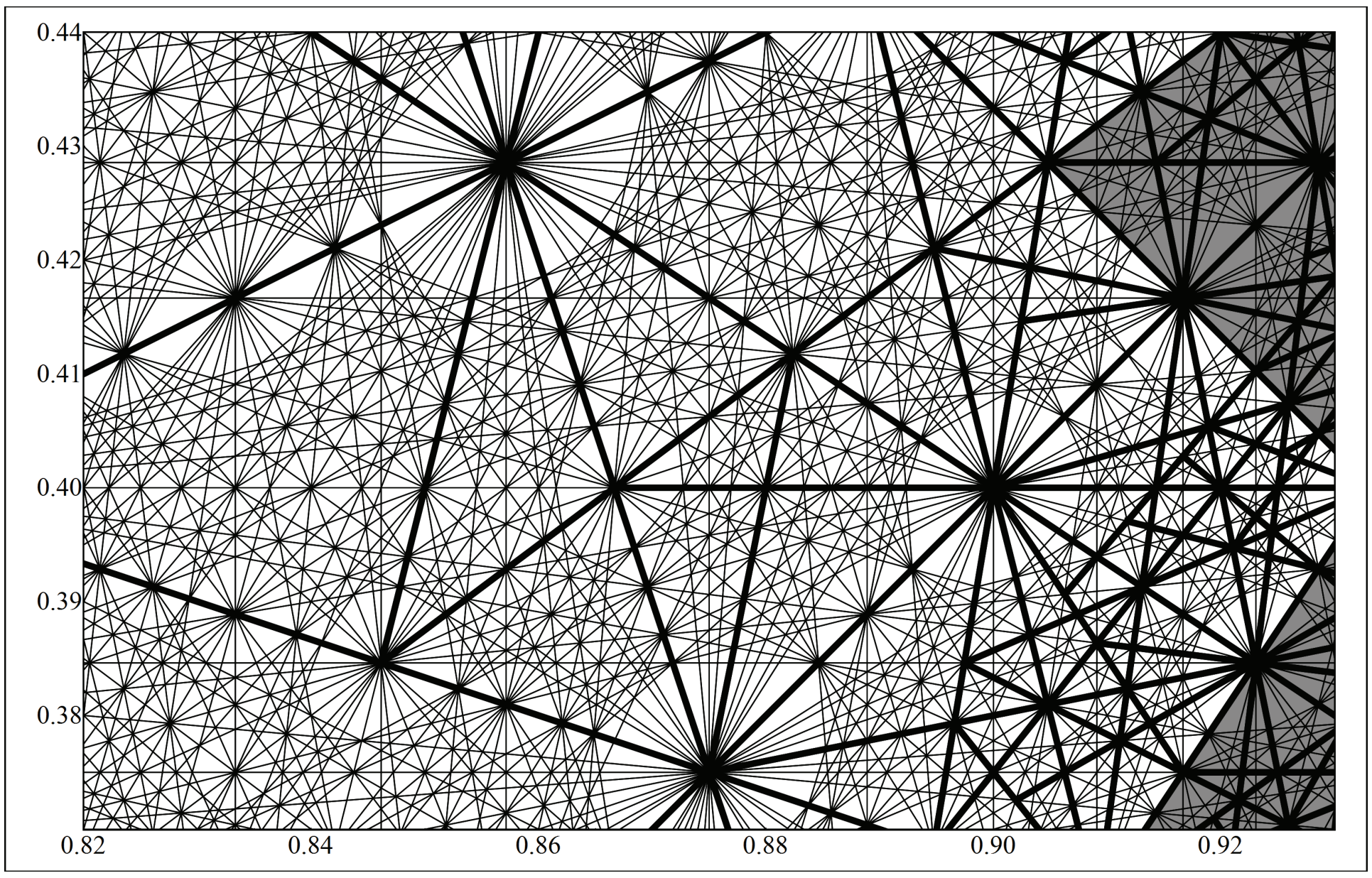

Fig. 1. Comparison of the decompositions obtained by Algorithm 1 (bold) and Algorithm 2 of $H = \mathrm{conv}\left(\left\{\left(\frac{41}{50}, \frac{37}{100}\right), \left(\frac{93}{100}, \frac{37}{100}\right), \left(\frac{93}{100}, \frac{11}{25}\right), \left(\frac{41}{50}, \frac{11}{25}\right)\right\}\right)$. Dark regions are contained in cutout polygons.

## 4. Results

Using Algorithm 2 of Section 3 the list of Appendix A has been found. It completely characterizes the set $\mathcal{D}_2^{(0)} \cap C$ where

$$\begin{aligned} C &:= C_1 \setminus C_2, \\ C_1 &:= \left\{(x,y) \in \mathbb{R}^2 \mid x \leq 1 - L\right\}, \\ C_2 &:= \operatorname{int}\Big(\operatorname{conv}\Big(\Big\{\Big(1-K, 2-K\Big), \Big(1-K+\sqrt{2}L, 2-K\Big), \\ &\qquad \Big(1-\sqrt{2}L, 2-2\sqrt{2}L\Big), \Big(1,2\Big)\Big\}\Big)\Big) \qquad (4.1) \end{aligned}$$

and $K = \frac{1}{20}$, $L = \frac{1}{512}$, which is considerably larger than what has been achieved with Brunotte's algorithm for regions so far ($L = \frac{1}{100}$, [18]). Note that $C_2$ is a small open quadrangle of width $L$ touching the boundary of $\mathcal{D}_2$ left of $(1,2)$. The convex sets used as inputs for the algorithm were the closed squares $\left[\frac{x}{n}, \frac{x+1}{n}\right] \times \left[\frac{y}{n}, \frac{y+1}{n}\right]$, where $n = 8192$, $(x,y) \in \mathbb{Z}^2$ with $\left\lfloor \frac{2n}{3} \right\rfloor \leq x \leq L - \frac{1}{n}$ and $-\frac{n}{2} \leq y \leq \frac{3n}{2} - 1$. The remaining regions have already been characterized in [2] (especially by Theorem 4.8 there, which covers the region $\left\{(x,y) \in \mathbb{R}^2 \mid 0 < x < 1 \wedge 0 < y < x+1 \wedge 4x < y^2 \wedge y > \frac{x}{\gamma_6} + \gamma_6\right\}$ where $\gamma_q$ is the positive root of $qt^3 + qt^2 - qt - q + 1$, $q \in \mathbb{N}$, and therefore reaches the boundary of $\mathcal{D}_2$) or were also treated by Algorithm 2 using suitable convex input sets. Every 5-tuple $(n, x, y, a_1, a_2)$ in the list of Appendix A represents a cutout polygon $P$ in the following way: Let $\mathbf{r} := (\frac{x}{n}, \frac{y}{n})$, $\mathbf{a} := (a_1, a_2)$, $m := \min\left\{k \in \mathbb{N} \mid \tau_{\mathbf{r}}^k(\mathbf{a}) = \mathbf{a}\right\}$, and $\pi := (\tau_{\mathbf{r}}(\mathbf{a}), \ldots, \tau_{\mathbf{r}}^m(\mathbf{a}))$. Then $P := P(\pi)$.

If $\{P_1, \ldots, P_{598}\}$ is the set of the 598 cutout polygons then

$$\mathcal{D}_2^{(0)} \cap C = \left\{(x,y) \in \mathbb{R}^2 \mid x \leq 1 \ \wedge\ |y| \leq x+1\right\} \cap C \setminus \bigcup_{k=1}^{598} P_k. \qquad (4.2)$$

Note that $\left\{(x,y) \in \mathbb{R}^2 \mid x \leq 1 \ \wedge\ |y| \leq x+1\right\}$ is the topological closure of $\mathcal{D}_2$ [1]. Also note that none of the given cutout polygons is redundant, which can easily be verified as any of the given parameters is contained solely in the corresponding cutout polygon.

The analysis of the list of cutout polygons leads to the following

**Theorem 4.1.**

*(i)* $\mathcal{D}_2^{(0)}$ *has at least 22 connected components*

*(ii) The largest connected component of* $\mathcal{D}_2^{(0)}$ *has at least 3 holes*

**Proof.** The parameters

$(\frac{1}{2}, \frac{1}{2})$, $(\frac{152}{157}, \frac{193}{157})$, $(\frac{313}{315}, \frac{239}{210})$, $(\frac{167}{168}, \frac{255}{224})$, $(\frac{314}{317}, \frac{359}{317})$, $(\frac{453}{455}, \frac{496}{455})$, $(\frac{305}{306}, \frac{37}{34})$, $(\frac{362}{363}, \frac{259}{242})$, $(\frac{356}{357}, \frac{382}{357})$, $(\frac{358}{359}, \frac{384}{359})$, $(\frac{1121}{1124}, \frac{601}{562})$, $(\frac{1375}{1378}, \frac{640}{689})$, $(\frac{2061}{2066}, \frac{959}{1033})$, $(\frac{309}{310}, \frac{141}{155})$, $(\frac{1533}{1538}, \frac{699}{769})$, $(\frac{989}{992}, \frac{901}{992})$, $(\frac{1127}{1133}, \frac{1009}{1133})$, $(\frac{1607}{1612}, \frac{691}{806})$, $(\frac{694}{697}, \frac{521}{697})$, $(\frac{92}{93}, \frac{16}{31})$, $(\frac{537}{539}, \frac{67}{539})$, $(\frac{304}{305}, \frac{38}{305})$

are contained in 22 distinct connected components of $\mathcal{D}_2^{(0)}$.

The parameters

$$\left(\tfrac{911}{914}, \tfrac{391}{457}\right), \left(\tfrac{2455}{2463}, \tfrac{2108}{2463}\right), \left(\tfrac{265}{266}, \tfrac{1}{4}\right)$$

are contained in 3 distinct holes of the largest connected component of $\mathcal{D}_2^{(0)}$. □

The figures below show the calculated cutout polygons and the resulting shape of $\mathcal{D}_2^{(0)}$ in the corresponding region. Figure 2 gives an overview of $\mathcal{D}_2^{(0)}$ and shows the regions which lie above and below the point $(1, 1)$, Fig. 3 shows several connected components, and Fig. 4 gives an example of holes.

The cutout polygons are represented in the following way: If an edge belongs to the polygon it is plotted solid, and dotted otherwise. Belonging of a vertex is indicated by a prominent dot at the respective position. In the images which show the resulting shape of $\mathcal{D}_2^{(0)}$ black regions do belong to $\mathcal{D}_2^{(0)}$, white regions do not and gray regions are not settled by now.

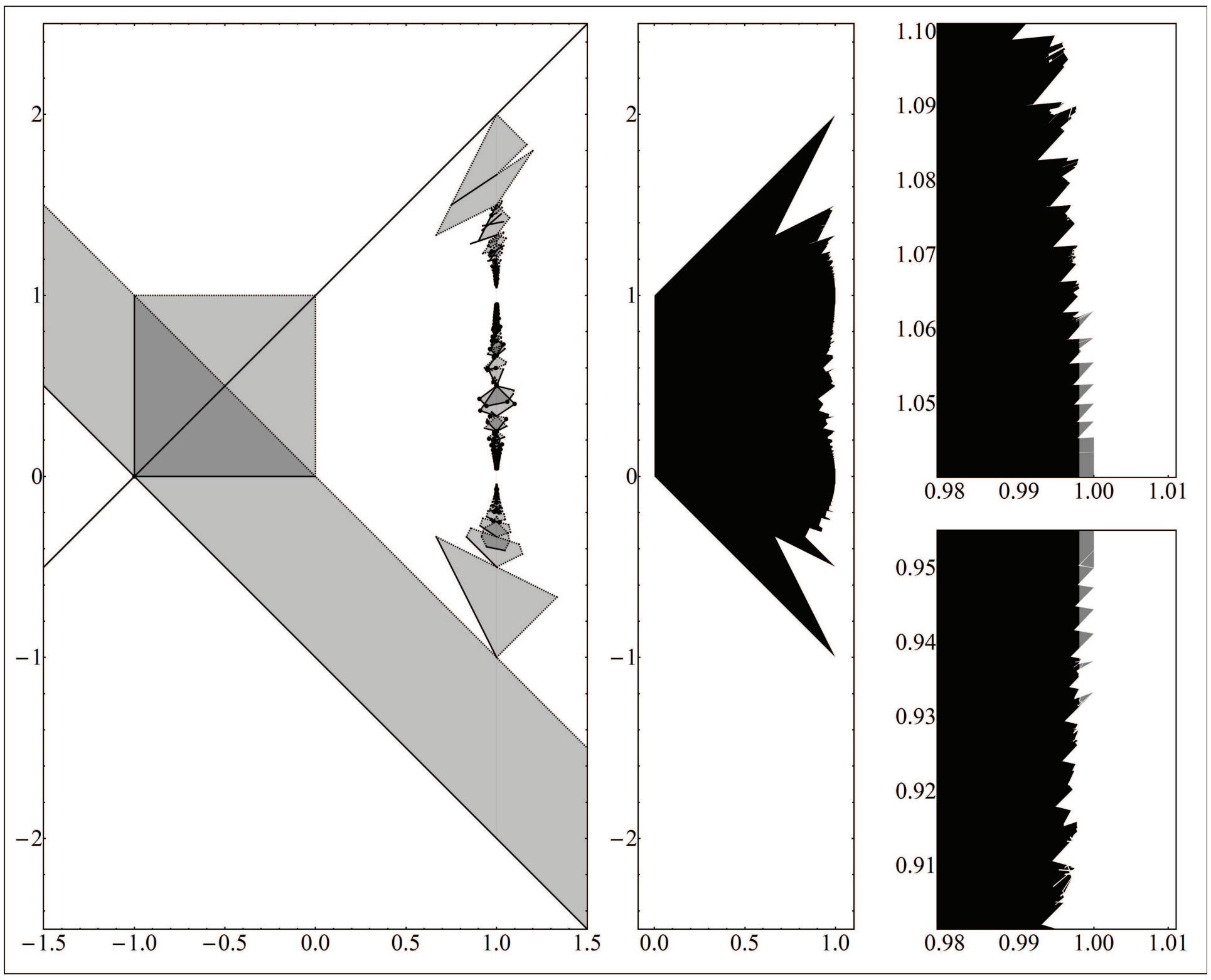

Fig. 2. Overview of $\mathcal{D}_2^{(0)}$ and the regions above and below $(1, 1)$.

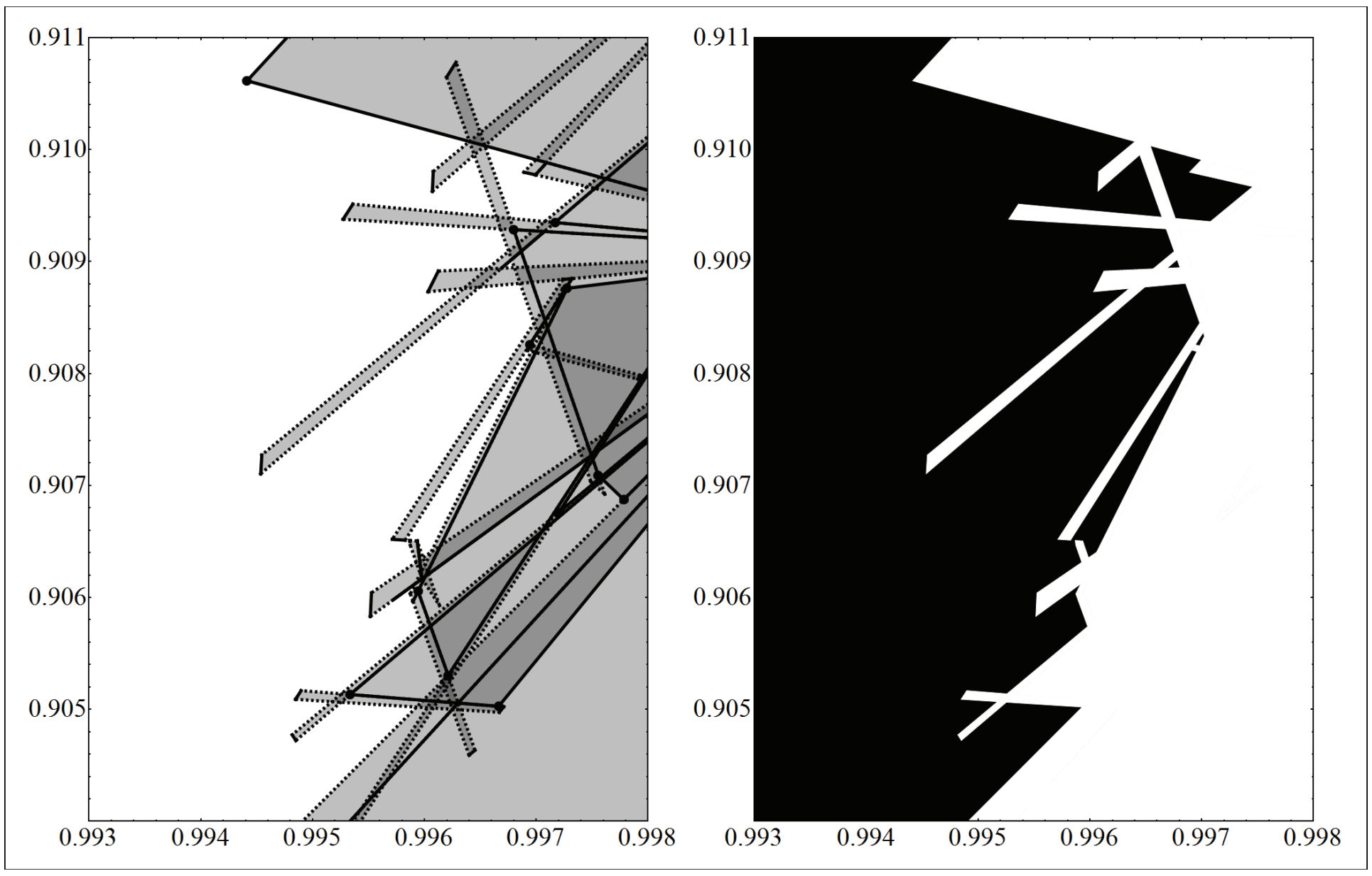

Fig. 3. Four connected components of $\mathcal{D}_2^{(0)}$.

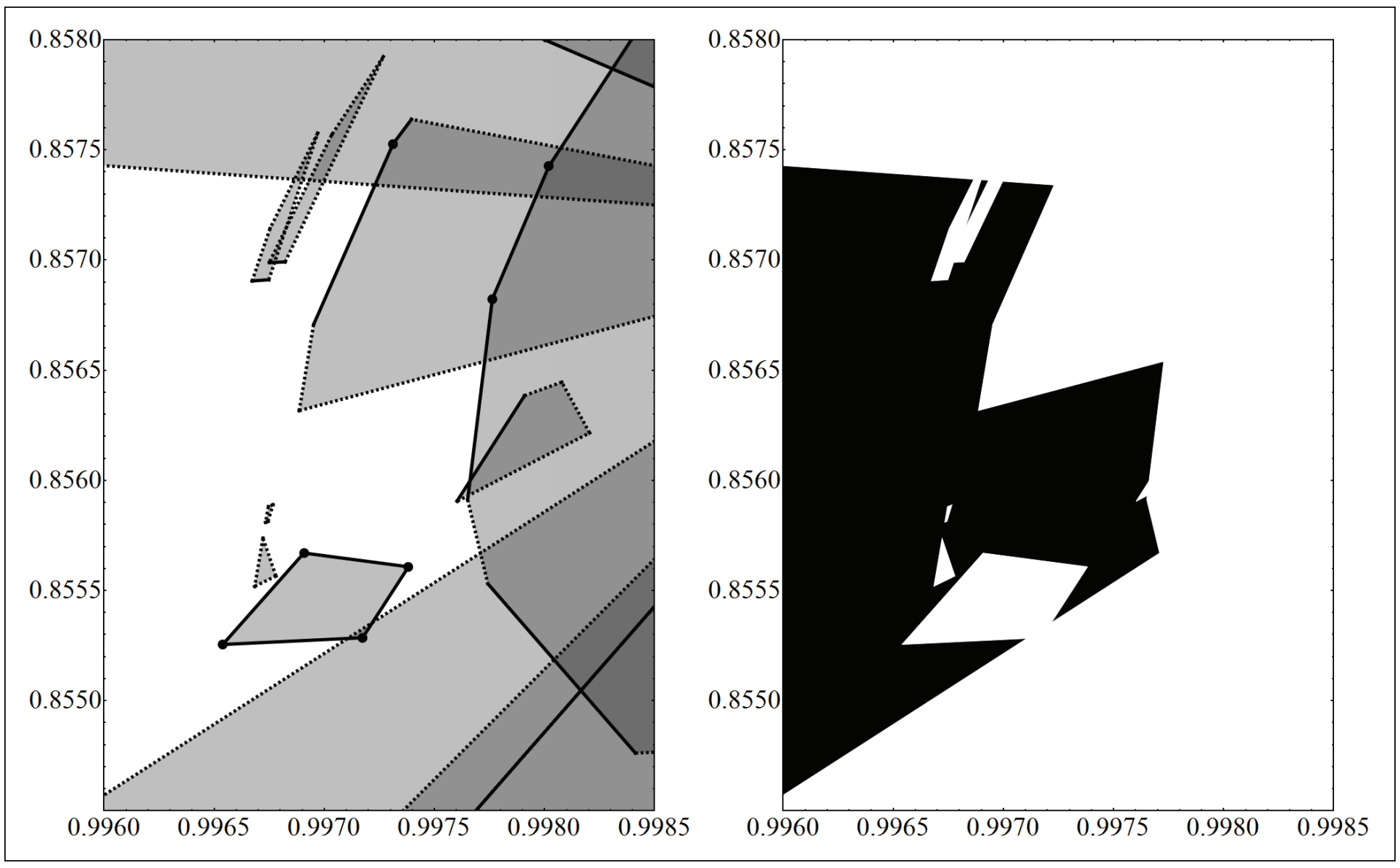

Fig. 4. Two holes of $\mathcal{D}_2^{(0)}$.

## 5. Infinite Families of Cutout Polygons

In [18] two infinite families of cutout polygons of $\mathcal{D}_2^{(0)}$ are deduced. The following Lemma provides a formal method to identify families of cutout polyhedra. It was used to find four additional infinite families of cutout polygons of $\mathcal{D}_2^{(0)}$ which - together with the two families already known and two singular cutout polygons - form a chain leading from $(1,-1)$ to $(1,2)$ along the "critical line" $\{(x,y) \in \mathbb{R}^2 \mid x = 1\}$. They cut out a considerable part of the region in $\mathcal{D}_2$ which has not been settled algorithmically in Section 4.

**Lemma 5.1.** *Let $\mathcal{H}$ denote a finite set of half-spaces in $\mathbb{R}^d$ and $P \in \mathcal{P}_d$ be bounded. Furthermore let $\mathcal{H}_o := \{H \in \mathcal{H} \mid H \text{ open}\}$, $\mathcal{H}_c := \{H \in \mathcal{H} \mid H \text{ closed}\}$, $\mathcal{F}_o(P) := \mathcal{F}^\circ(P) \setminus \mathcal{P}(P)$, $\mathcal{F}_c(P) := \mathcal{F}^\circ(P) \cap \mathcal{P}(P)$ and $A_P := \bigcap \{A \text{ affine subspace of } \mathbb{R}^d \mid P \subseteq A\}$. Then $P = \bigcap \mathcal{H}$ iff the following holds:*

*(i)* $\forall\, F \in \mathcal{F}^\circ(P) : |F| = 1 \Rightarrow \forall\, H \in \mathcal{H} : F \subseteq \overline{H}$
*(ii)* $A_P = \mathbb{R}^d \vee \exists\, \mathcal{H}' \subseteq \mathcal{H} : A_P = \bigcap \mathcal{H}'$
*(iii)* $\forall\, F \in \mathcal{F}_o(P) : \exists\, H \in \mathcal{H}_o : F \subseteq \partial H$
*(iv)* $\forall\, F \in \mathcal{F}_c(P) : \overline{F} \neq \overline{P} \Rightarrow \exists\, H \in \mathcal{H}_c : F \subseteq \partial H \wedge P \nsubseteq \partial H$
*(v)* $\forall\, F \in \mathcal{F}_c(P) : \nexists\, H \in \mathcal{H}_o : F \subseteq \partial H$

**Proof.** It is obvious that $P = \bigcap \mathcal{H} \Rightarrow$ (i)$\wedge \ldots \wedge$(v). In the other direction the boundedness of $P$, (i) and (v) imply that $P \subseteq \bigcap \mathcal{H}$. Also it follows from (ii) that $\bigcap \mathcal{H} \subseteq A_P$ – the smallest affine subspace containing $P$ (if $P$ is nonempty, otherwise the statement of the Lemma is trivial due to (ii)). In addition to the restriction $P \nsubseteq \partial H$ in (iv) this allows to assume w.l.o.g. that $P$ is nondegenerate ($\Leftrightarrow A_P = \mathbb{R}^d$). But then it follows from (iii) and (iv) that $\bigcap \mathcal{H} \subseteq P$. □

**Definition 5.2.** *For two finite tuples $S$ and $T$ let $S \sqcup T$ denote the tuple obtained by concatenation of $S$ and $T$.*

*For a tuple $(T_1, \ldots, T_n)$ of $n \in \mathbb{N}$ finite tuples let* $\operatorname{shuffle}(T_1, \ldots, T_n)$ *denote the tuple obtained by successively stringing together the first entries of the tuples (in the given order) followed by the second entries and so forth, with tuples having too little entries being skipped (e.g.* $\operatorname{shuffle}((1,2),(3),(4,5,6)) = (1,3,4,2,5,6)$*).*

Using the notions of Definition 5.2 one can define the following families of cycles (for convenience the families studied in [18] are also given – $C_2(n)$ and $C_6(n)$):

$$\begin{aligned} C_0(1) &:= ((-3,3),(3,-2),(-2,1),(1,1),(1,-2),(-2,3),(3,-3)) \\ C_0(2) &:= ((-5,1),(1,5),(5,-3),(-3,-3),(-3,5),(5,1),(1,-5),(-5,2),(2,4), \\ &\quad\ (4,-4),(-4,-1),(-1,5),(5,-1),(-1,-4),(-4,4),(4,2),(2,-5)) \end{aligned}$$

$$
\begin{aligned}
C_1^{(1)}(n) &:= ((-2n,2k))_{k=1}^{n} \sqcup ((-2n+2k,2n))_{k=1}^{n-1} \sqcup ((2k-1,2n-2k))_{k=1}^{n-1} \sqcup \\
&\quad ((2n-1,-2k+1))_{k=1}^{n} \sqcup ((2n-2k-1,-2n+1))_{k=1}^{n-1} \sqcup \\
&\quad ((-2k,-2n+2k+1))_{k=1}^{n-1} \\
C_1^{(2)}(n) &:= ((2k,2n-2k))_{k=1}^{n-1} \sqcup ((2n,-2k+1))_{k=1}^{n} \sqcup ((2n-2k,-2n+1))_{k=1}^{n-1} \sqcup \\
&\quad ((-2k+1,-2n+2k+1))_{k=1}^{n-1} \sqcup ((-2n+1,2k))_{k=1}^{n} \sqcup \\
&\quad ((-2n+2k+1,2n))_{k=1}^{n-1} \\
C_1^{(3)}(n) &:= ((2n-2k,-2n))_{k=1}^{n-1} \sqcup ((-2k+1,-2n+2k))_{k=1}^{n-1} \sqcup \\
&\quad ((-2n+1,2k-1))_{k=1}^{n} \sqcup ((-2n+2k+1,2n-1))_{k=1}^{n-1} \sqcup \\
&\quad ((2k,2n-2k-1))_{k=1}^{n-1} \sqcup ((2n,-2k))_{k=1}^{n} \\
C_1(n) &:= \text{shuffle}(C_1^{(1)}(n), C_1^{(2)}(n), C_1^{(3)}(n)),\ n \geq 2
\end{aligned}
$$

$$
\begin{aligned}
C_2^{(1)}(n) &:= ((-2n,2k-1))_{k=1}^{n+1} \sqcup ((-2n+2k,2n+1))_{k=1}^{n-1} \\
C_2^{(2)}(n) &:= ((2k-1,2n-2k+1))_{k=1}^{n} \sqcup ((2n+1,-2k))_{k=1}^{n} \\
C_2^{(3)}(n) &:= ((2n-2k+1,-2n))_{k=1}^{n} \sqcup ((-2k,-2n+2k))_{k=1}^{n-1} \\
C_2(n) &:= \text{shuffle}(C_2^{(1)}(n), C_2^{(2)}(n), C_2^{(3)}(n)),\ n \geq 1
\end{aligned}
$$

$$
\begin{aligned}
C_3^{(1)}(n) &:= ((-2n-1,1)) \sqcup ((-2n+2k-2,-2k))_{k=1}^{n} \sqcup ((2k-1,-2n-1))_{k=1}^{n} \\
C_3^{(2)}(n) &:= ((1,2n+1)) \sqcup ((-2k,2n+2))_{k=1}^{n-1} \sqcup ((-2n,2n+1)) \sqcup \\
&\quad ((-2n-1,2n-2k+1))_{k=1}^{n-1} \\
C_3^{(3)}(n) &:= ((2n+1,-2n)) \sqcup ((2n+2,-2n+2k))_{k=1}^{n-1} \sqcup ((2n-2k+3,2k-1))_{k=1}^{n} \\
C_3(n) &:= \text{shuffle}(C_3^{(1)}(n), C_3^{(2)}(n), C_3^{(3)}(n)),\ n \geq 2
\end{aligned}
$$

$$
\begin{aligned}
C_4^{(1)}(n) &:= ((-2n-1,2)) \sqcup ((-2n+2k-2,-2k+1))_{k=1}^{n} \sqcup ((2k-1,-2n))_{k=1}^{n-1} \sqcup \\
&\quad ((2n-1,-2n+1)) \sqcup ((2n,-2n+2k+1))_{k=1}^{n-1} \sqcup ((2n-2k+1,2k))_{k=1}^{n} \sqcup \\
&\quad ((-2k,2n+1))_{k=1}^{n-1} \sqcup ((-2n,2n)) \sqcup ((-2n-1,2n-2k))_{k=1}^{n-2} \\
C_4^{(2)}(n) &:= ((2,2n)) \sqcup ((-2k+1,2n+1))_{k=1}^{n-1} \sqcup ((-2n+1,2n)) \sqcup \\
&\quad ((-2n,2n-2k))_{k=1}^{n-1} \sqcup ((-2n+2k-1,-2k+1))_{k=1}^{n} \sqcup ((2k,-2n))_{k=1}^{n-1} \sqcup \\
&\quad ((2n,-2n+1)) \sqcup ((2n+1,-2n+2k+1))_{k=1}^{n-1} \sqcup ((2n-2k+2,2k))_{k=1}^{n-1} \\
C_4^{(3)}(n) &:= ((2n,-2n)) \sqcup ((2n+1,-2n+2k))_{k=1}^{n-1} \sqcup ((2n-2k+2,2k-1))_{k=1}^{n} \sqcup \\
&\quad ((-2k+1,2n))_{k=1}^{n-1} \sqcup ((-2n+1,2n-1)) \sqcup ((-2n,2n-2k-1))_{k=1}^{n-1} \sqcup \\
&\quad ((-2n+2k-1,-2k))_{k=1}^{n} \sqcup ((2k,-2n-1))_{k=1}^{n-1} \\
C_4(n) &:= \text{shuffle}(C_4^{(1)}(n), C_4^{(2)}(n), C_4^{(3)}(n)),\ n \geq 2
\end{aligned}
$$

$$
\begin{aligned}
C_5^{(1)}(n) &:= ((-n-1,1)) \sqcup ((-n+k-1,k+2))_{k=1}^{n-1} \\
C_5^{(2)}(n) &:= ((1,n+1)) \sqcup ((k+2,n-k+1))_{k=1}^{n-2} \sqcup ((n+1,1)) \\
C_5^{(3)}(n) &:= ((n-k+2,-k-1))_{k=1}^{n-1} \sqcup ((1,-n-1)) \\
C_5^{(4)}(n) &:= ((-k-1,-n+k-1))_{k=1}^{n-1} \\
C_5(n) &:= \text{shuffle}(C_5^{(1)}(n), C_5^{(2)}(n), C_5^{(3)}(n), C_5^{(4)}(n)),\ n \geq 2
\end{aligned}
$$

$$
\begin{aligned}
C_6^{(1)}(n) &:= ((-n+k-1,-k))_{k=1}^{n} \sqcup ((1,-n))\\
C_6^{(2)}(n) &:= ((-k,n-k+1))_{k=1}^{n} \sqcup ((n-k+1,k+1))_{k=1}^{n}\\
C_6^{(3)}(n) &:= ((k+1,-n+k))_{k=1}^{n-1} \sqcup ((n+1,1))\\
C_6(n) &:= \mathrm{shuffle}(C_6^{(1)}(n), C_6^{(2)}(n), C_6^{(3)}(n)),\ n \geq 1
\end{aligned}
$$

Using Lemma 5.1 one can verify that the corresponding cutout polygons of the cycles defined above have the shapes given in Proposition 5.3 below. Every polygon is given by a list of its vertices in counterclockwise order. A vertex is overlined iff it belongs to the respective polygon. Belonging of an edge is indicated by a solid (contained) or dotted (not contained) line between the endpoints.

**Proposition 5.3.** *With the notions defined above we have the following infinite families of cutout polygons:*

$C_0(n)$, $n=1$: $\left(\frac{3}{4},\frac{3}{2}\right)$—$\left(1,\frac{5}{3}\right)$- -$\left(\frac{7}{6},\frac{11}{6}\right)$- -$\left(1,2\right)$- -

$n=2$: $\left(\frac{25}{26},\frac{15}{26}\right)$—$\left(1,\frac{1}{2}\right)$—$\left(\frac{28}{27},\frac{16}{27}\right)$- -$\left(1,\frac{3}{5}\right)$- -

$C_1(n)$, $n\geq 2$: $\left(1-\frac{1}{4n^2-4n+2},1+\frac{2n-1}{4n^2-4n+2}\right)$- -$\left(1,1+\frac{1}{2n-1}\right)$- -
$\left(1+\frac{1}{4n^2-2},1+\frac{2n+2}{4n^2-2}\right)$- -$\left(1,1+\frac{1}{2n-2}\right)$- -

$C_2(n)$, $n=1$: $\left(\frac{2}{3},\frac{4}{3}\right)$- -$\left(1,\frac{3}{2}\right)$- -$\left(\frac{6}{5},\frac{9}{5}\right)$- -$\left(\frac{3}{4},\frac{3}{2}\right)$- -

$n\geq 2$: $\left(1-\frac{1}{4n^2-2n+1},1+\frac{2n-1}{4n^2-2n+1}\right)$- -$\left(1,1+\frac{1}{2n}\right)$- -
$\left(1+\frac{1}{4n^2+2n-1},1+\frac{2n+2}{4n^2+2n-1}\right)$- -$\left(1,1+\frac{1}{2n-1}\right)$- -

$C_3(n)$, $n\geq 2$: $\overline{\left(1-\frac{1}{4n^2+6n-1},1-\frac{2n+4}{4n^2+6n-1}\right)}$—$\overline{\left(1,1-\frac{1}{2n-1}\right)}$—
$\overline{\left(1+\frac{1}{4n^2+6n-2},1-\frac{2n+3}{4n^2+6n-2}\right)}$—$\overline{\left(1,1-\frac{1}{2n}\right)}$—

$C_4(n)$, $n=2$: $\overline{\left(\frac{19}{20},\frac{3}{5}\right)}$—$\overline{\left(\frac{21}{22},\frac{13}{22}\right)}$—$\left(1,\frac{3}{5}\right)$- -$\left(\frac{22}{21},\frac{13}{21}\right)$- -$\left(\frac{20}{19},\frac{12}{19}\right)$- -$\left(1,\frac{2}{3}\right)$—

$n\geq 3$: $\overline{\left(1-\frac{1}{4n^2+4n-4},1-\frac{2n+4}{4n^2+4n-4}\right)}$—$\left(1,1-\frac{1}{2n-2}\right)$- -
$\left(1+\frac{1}{4n^2+4n-5},1-\frac{2n+3}{4n^2+4n-5}\right)$- -$\left(1,1-\frac{1}{2n-1}\right)$—

$C_5(n)$, $n=2$: $\overline{\left(\frac{10}{11},\frac{4}{11}\right)}$—$\left(1,\frac{1}{3}\right)$—$\overline{\left(\frac{11}{10},\frac{2}{5}\right)}$—$\overline{\left(1,\frac{1}{2}\right)}$—

$n=3$: $\left(\frac{14}{15},\frac{4}{15}\right)$—$\left(1,\frac{1}{4}\right)$—$\left(\frac{19}{18},\frac{5}{18}\right)$- -$\left(1,\frac{1}{3}\right)$- -

$n=4$: $\left(\frac{22}{23},\frac{5}{23}\right)$—$\overline{\left(\frac{23}{24},\frac{5}{24}\right)}$—$\left(1,\frac{1}{5}\right)$—$\left(\frac{24}{23},\frac{5}{23}\right)$- -$\left(1,\frac{1}{4}\right)$- -

$n\geq 5$: $\left(1-\frac{1}{n^2+n+3},\frac{n+1}{n^2+n+3}\right)$—$\overline{\left(1-\frac{1}{n^2+2n},\frac{n+1}{n^2+2n}\right)}$—$\left(1,\frac{1}{n+1}\right)$—
$\overline{\left(1+\frac{1}{n^2+2n-1},\frac{n+1}{n^2+2n-1}\right)}$—$\left(1+\frac{1}{n^2+n+3},\frac{n+1}{n^2+n+3}\right)$- -$\left(1,\frac{1}{n}\right)$- -

$C_6(n)$, $n=1$ $\left(\frac{2}{3},-\frac{1}{3}\right)$—$\left(1,-1\right)$- -$\left(\frac{4}{3},-\frac{2}{3}\right)$- -

$n\geq 2$: $\left(1-\frac{1}{n^2+2},-\frac{n}{n^2+2}\right)$—$\left(1,-\frac{1}{n}\right)$- -$\left(1+\frac{1}{n^2+n+1},-\frac{n+1}{n^2+n+1}\right)$- -
$\left(1+\frac{1}{n^2+2n},-\frac{n+1}{n^2+2n}\right)$- -$\left(1-\frac{1}{n^2+n+1},-\frac{n}{n^2+n+1}\right)$- -

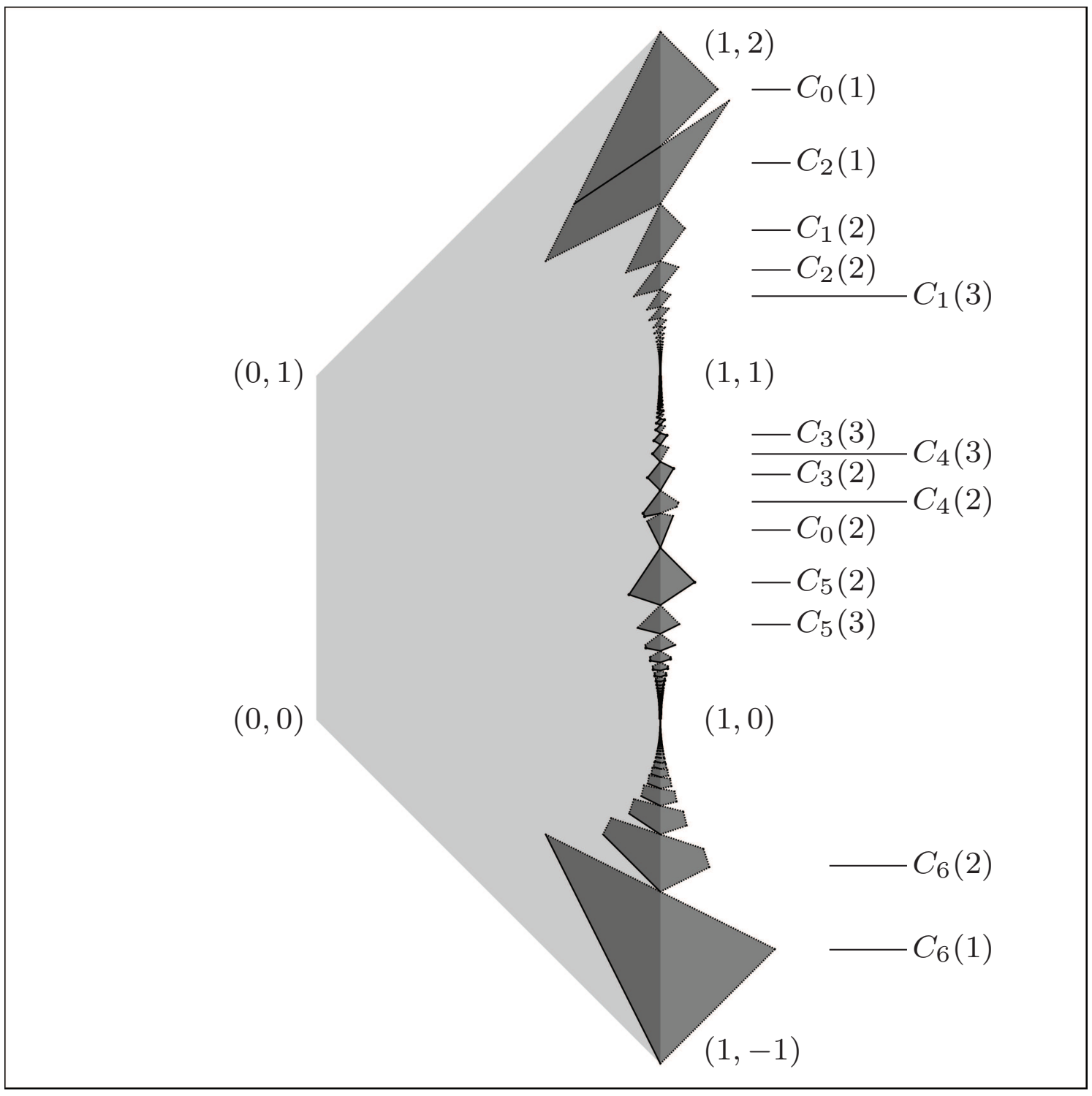

Fig. 5. Six families of cutout polygons.

## 6. Further Remarks and Conjectures

The answers to the topological questions given in this paper of course raise new questions. How many connected components and holes does $\mathcal{D}_2^{(0)}$ have? Considering its highly complicated structure, the answer might be "infinitely many" in both cases. A full characterization of $\mathcal{D}_2^{(0)}$ seems to be out of reach for the moment, but the two "critical points" (cf. [1]) $(1, 0)$ and $(1, 1)$ appear to be asymmetric in terms of complexity of the surrounding areas. While the latter hoards several connected components and holes, the former seems to have a comparably simple, maybe even regular neighborhood. It might be possible to characterize $\mathcal{D}_2^{(0)}$ by finitely many classes of cutout polygons at least for $y \leq 1 - \epsilon$ where $\epsilon > 0$. Furthermore, besides of some singular results (e.g. [10,11]), the structure of $\mathcal{D}_d^{(0)}$ for $d > 2$ is still open.

## Acknowledgments

The author is supported by the Austrian Science Fund (FWF): W1230, Doctoral Program "Discrete Mathematics".

The author wishes to dedicate this paper to his newborn daughter Mia Kallisto!

## Appendix A. List of Cutout Polygons

(1, 0, -1, 1, 1) (1, 0, 1, -1, 1) (3, -1, 1, 0, 1) (3, 3, -2, 2, 1) (3, 3, -1, 5, 1) (3, 3, 2, 5, 1) (4, 4, 1, 6, 1) (4, 4, 7, 1, 1)
(5, 5, 1, 11, 4) (5, 5, 8, 1, 1) (6, 6, -1, 9, 1) (6, 6, 1, 10, 1) (6, 6, 5, 7, 1) (7, 7, -3, 3, 1) (7, 7, -1, 16, 4) (7, 7, 4, 5, 1)
(7, 7, 9, 3, 1) (8, 8, -1, 21, 4) (8, 8, 3, 3, 1) (8, 8, 5, 4, 1) (8, 8, 7, 9, 1) (8, 8, 9, 25, 2) (8, 8, 11, 2, 1) (9, 9, -2, 5, 1)
(9, 9, -1, 13, 1) (9, 9, 1, 14, 1) (9, 9, 7, 6, 1) (9, 9, 8, 11, 1) (9, 9, 10, 39, 4) (9, 9, 11, 4, 1) (10, 9, 13, 2, 1) (10, 10, -1, 29, 4)
(10, 10, 11, 15, 2) (11, 11, 10, 13, 1) (11, 11, 12, 27, 2) (11, 11, 13, 5, 1) (12, 11, 5, 3, 3) (12, 11, 16, 2, 1)
(12, 12, -1, 17, 1) (12, 12, 1, 18, 1) (13, 13, -1, 35, 4) (13, 13, 3, 5, 1) (13, 13, 11, 8, 1) (13, 13, 15, 6, 1)
(14, 13, -4, 4, 1) (14, 14, -1, 54, 2) (14, 14, 15, 30, 2) (15, 15, -2, 8, 1) (15, 15, -1, 21, 1) (15, 15, 1, 22, 1)
(15, 15, 2, 8, 1) (15, 15, 14, 17, 1) (15, 15, 16, 35, 2) (15, 15, 17, 7, 1) (16, 16, 17, 25, 2) (17, 17, 1, 26, 1)
(17, 17, 16, 19, 1) (17, 17, 19, 8, 1) (18, 18, -1, 25, 1) (19, 18, 5, 4, 1) (19, 19, 18, 21, 1) (19, 19, 21, 9, 1)
(21, 21, 19, 12, 1) (21, 21, 23, 10, 1) (22, 21, 28, 2, 4) (23, 23, 25, 11, 1) (25, 24, 17, 6, 1) (25, 25, 27, 12, 1)
(26, 25, 31, 4, 1) (27, 26, 9, 9, 1) (27, 26, 35, 8, 2) (27, 27, -2, 14, 1) (27, 27, 25, 15, 1) (27, 27, 29, 13, 1)
(28, 27, 34, 7, 2) (29, 29, 27, 16, 1) (29, 29, 31, 14, 1) (30, 29, 36, 4, 1) (30, 29, 37, 7, 1) (31, 29, 43, 4, 1)
(31, 31, 33, 15, 1) (32, 30, 19, 3, 2) (32, 31, -6, 6, 1) (32, 32, -3, 11, 1) (33, 33, -2, 17, 1) (33, 33, 2, 17, 1)
(33, 33, 31, 35, 2) (33, 33, 35, 16, 1) (34, 33, 6, 6, 1) (34, 33, 42, 7, 2) (35, 34, 43, 7, 2) (35, 35, 37, 17, 1)
(37, 35, 22, 4, 2) (37, 37, 35, 20, 1) (37, 37, 39, 18, 1) (38, 37, 56, 15, 3) (39, 39, -2, 20, 1) (39, 39, 2, 20, 1)
(39, 39, 41, 19, 1) (40, 39, 29, 8, 1) (40, 39, 30, 10, 2) (41, 38, 56, 4, 2) (41, 41, -2, 21, 1) (41, 41, 2, 21, 1)
(41, 41, 43, 20, 1) (43, 43, -2, 22, 1) (43, 43, 2, 22, 1) (43, 43, 45, 21, 1) (44, 43, -7, 7, 1) (45, 45, -2, 23, 1)
(45, 45, 2, 23, 1) (46, 44, 27, 4, 2) (46, 45, 7, 7, 1) (49, 47, 60, 5, 2) (50, 49, 57, 6, 1) (50, 49, 58, 15, 2)
(52, 51, 27, 15, 1) (54, 53, 41, 9, 2) (55, 54, 11, 11, 1) (55, 54, 42, 10, 4) (55, 54, 64, 15, 2) (56, 55, 64, 6, 1)
(56, 55, 65, 15, 2) (61, 60, -15, 16, 2) (61, 60, 71, 15, 3) (61, 60, 72, 7, 4) (63, 62, 75, 9, 6) (67, 66, 55, 18, 3)
(68, 67, 45, 21, 2) (68, 67, 54, 13, 2) (68, 67, 56, 20, 4) (68, 67, 79, 14, 3) (68, 68, 5, 14, 1) (71, 69, -17, 9, 2)
(71, 70, 59, 11, 1) (72, 71, 83, 9, 3) (74, 72, 57, 5, 2) (74, 73, -9, 9, 1) (74, 73, 18, 17, 1) (74, 73, 59, 14, 2)
(74, 73, 86, 14, 2) (76, 74, 57, 9, 3) (76, 75, 9, 9, 1) (79, 76, 97, 6, 1) (79, 77, 59, 10, 3) (80, 78, 59, 8, 3)
(81, 80, 92, 10, 2) (82, 81, 66, 11, 3) (82, 81, 91, 8, 1) (83, 81, 64, 6, 2) (86, 85, 69, 13, 4) (89, 87, 48, 9, 2)
(90, 89, 100, 8, 1) (92, 91, -10, 10, 1) (94, 93, 10, 10, 1) (98, 97, 14, 14, 1) (98, 97, 111, 24, 2) (99, 97, 75, 8, 3)
(100, 99, 85, 13, 1) (114, 113, 11, 11, 1) (116, 115, 102, 10, 1) (120, 119, 136, 22, 3) (121, 120, 96, 21, 3)
(122, 121, 133, 10, 1) (122, 121, 137, 14, 1) (124, 123, 141, 16, 2) (127, 123, 183, 7, 1) (127, 126, 109, 14, 2)
(128, 127, 16, 15, 3) (129, 127, 106, 17, 4) (132, 130, 111, 7, 2) (132, 131, 144, 10, 1) (134, 133, -12, 12, 1)
(134, 133, 152, 26, 4) (136, 135, 12, 12, 1) (136, 135, 152, 13, 3) (137, 134, 74, 9, 2) (137, 136, 154, 14, 2)
(139, 138, 121, 16, 2) (140, 139, 155, 38, 3) (141, 139, 34, 19, 2) (141, 140, 20, 21, 2) (141, 140, 160, 26, 4)
(142, 141, 157, 14, 3) (144, 142, 115, 14, 2) (145, 143, 122, 8, 2) (147, 146, 163, 14, 3) (148, 146, -29, 15, 2)
(148, 147, 172, 22, 4) (151, 150, 134, 16, 1) (153, 151, -30, 15, 4) (154, 152, 37, 18, 2) (155, 153, 127, 19, 5)
(156, 155, 22, 28, 2) (158, 155, 183, 13, 3) (158, 157, -13, 13, 1) (159, 157, 127, 14, 2) (160, 159, 13, 13, 1)
(168, 167, 183, 28, 2) (170, 169, 183, 12, 1) (172, 170, 195, 11, 2) (173, 169, 129, 10, 3) (177, 176, 194, 16, 2)
(179, 177, 147, 19, 4) (179, 178, 195, 28, 2) (180, 178, 93, 22, 1) (180, 179, 196, 27, 2) (182, 181, 136, 44, 3)
(182, 181, 196, 12, 1) (183, 182, 201, 16, 3) (185, 183, 152, 19, 4) (185, 183, 211, 10, 4) (189, 187, 214, 23, 2)
(191, 189, 157, 19, 3) (191, 190, 208, 27, 2) (191, 190, 212, 19, 1) (192, 191, 209, 29, 3) (193, 192, 173, 18, 1)
(195, 192, 160, 18, 4) (200, 199, 222, 19, 1) (201, 200, 219, 27, 3) (202, 199, 106, 15, 1) (202, 201, 179, 21, 3)
(203, 202, 25, 24, 4) (203, 202, 178, 36, 4) (206, 204, 181, 9, 2) (212, 211, -15, 15, 1) (214, 213, 15, 15, 1)
(218, 215, 114, 16, 4) (220, 216, 255, 14, 2) (220, 217, 115, 16, 4) (220, 219, 202, 14, 1) (221, 218, 181, 18, 4)
(222, 221, 242, 54, 4) (223, 221, 196, 10, 2) (226, 225, 241, 14, 1) (230, 229, 207, 19, 2) (231, 230, 259, 23, 7)
(232, 231, 210, 62, 6) (233, 232, 174, 53, 3) (238, 237, 267, 35, 7) (240, 239, 256, 14, 1) (241, 240, 262, 18, 3)
(242, 239, 275, 11, 1) (242, 241, -16, 16, 1) (242, 241, 263, 20, 2) (244, 243, 16, 16, 1) (244, 243, 146, 48, 2)
(244, 243, 268, 24, 5) (245, 244, 222, 20, 2) (249, 248, 269, 19, 3) (250, 249, 25, 25, 1) (254, 251, 288, 12, 1)
(255, 254, 191, 63, 3) (255, 254, 231, 60, 7) (256, 255, 23, 22, 3) (258, 257, 32, 68, 3) (258, 257, 283, 24, 1)
(260, 259, 238, 21, 1) (265, 264, 33, 67, 3) (265, 264, 66, 66, 1) (265, 264, 238, 26, 2) (265, 264, 285, 25, 2)
(268, 267, 294, 24, 1) (270, 268, 235, 15, 3) (272, 270, 301, 38, 4) (273, 271, 302, 39, 4) (273, 272, 34, 68, 3)
(275, 273, 239, 15, 2) (276, 275, 248, 27, 2) (277, 275, 241, 16, 3) (277, 276, 298, 25, 2) (278, 276, 243, 15, 3)
(278, 277, 250, 27, 2) (278, 277, 299, 25, 2) (278, 277, 301, 20, 3) (279, 277, 242, 29, 6) (280, 279, 301, 20, 2)
(281, 280, -28, 29, 2) (286, 285, -57, 56, 4) (290, 288, 321, 39, 3) (290, 289, 307, 16, 1) (293, 289, 240, 18, 4)
(295, 294, 321, 30, 2) (296, 294, 267, 11, 2) (299, 298, 262, 27, 3) (301, 300, 322, 95, 7) (304, 303, 260, 79, 5)
(306, 305, 324, 16, 1) (308, 306, 341, 39, 5) (308, 307, -18, 18, 1) (309, 308, 278, 26, 4) (310, 309, 18, 18, 1)
(317, 315, 286, 12, 2) (318, 317, 293, 28, 2) (327, 325, 284, 16, 3) (327, 326, 297, 31, 4) (328, 327, 302, 23, 1)
(329, 328, 299, 33, 3) (331, 330, 354, 96, 6) (335, 333, 291, 17, 3) (335, 334, 287, 45, 1) (337, 336, 303, 31, 5)
(338, 337, 386, 46, 1) (341, 336, 281, 20, 4) (343, 341, 298, 17, 2) (344, 343, -19, 19, 1) (346, 345, 19, 19, 1)
(351, 349, 389, 19, 1) (352, 351, 32, 64, 2) (352, 351, 402, 47, 1) (353, 352, 321, 33, 2) (353, 352, 411, 54, 4)
(354, 353, 379, 24, 1) (356, 355, 378, 41, 2) (356, 355, 381, 95, 4) (360, 359, 329, 48, 2) (361, 360, 386, 23, 3)
(362, 361, 381, 18, 1) (363, 361, 403, 35, 3) (363, 362, 335, 24, 2) (363, 362, 388, 49, 3) (365, 362, 317, 15, 2)
(365, 364, 331, 61, 6) (367, 365, 411, 21, 1) (368, 367, 394, 24, 2) (369, 368, 395, 30, 2) (370, 369, 396, 30, 2)
(372, 370, 413, 34, 5) (372, 371, 398, 95, 6) (373, 372, 396, 40, 2) (374, 372, 333, 44, 7) (374, 373, 347, 102, 9)
(377, 376, 403, 49, 3) (380, 379, 400, 18, 1) (381, 380, 406, 24, 3) (381, 380, 435, 53, 1) (385, 384, 410, 24, 3)
(386, 384, 55, 27, 5) (389, 387, 353, 32, 3) (389, 388, 413, 40, 2) (390, 388, 353, 62, 7) (390, 389, 414, 42, 3)
(391, 390, 423, 69, 4) (393, 392, 49, 67, 4) (395, 393, 350, 20, 4) (400, 399, 373, 26, 1) (401, 400, 50, 67, 4)
(402, 400, 369, 13, 2) (403, 401, 357, 21, 4) (403, 402, 31, 30, 4) (403, 402, 439, 34, 5) (404, 403, 429, 40, 3)
(404, 403, 437, 69, 3) (405, 404, 376, 104, 8) (406, 405, 431, 39, 3) (407, 404, 353, 29, 5) (408, 407, 378, 26, 2)
(409, 408, 29, 29, 1) (409, 408, 436, 25, 3) (411, 410, 381, 33, 2) (416, 415, 445, 95, 6) (416, 415, 450, 69, 3)
(417, 409, 483, 14, 2) (423, 422, 387, 33, 1) (423, 422, 461, 38, 1) (427, 425, 392, 14, 2) (428, 425, 380, 43, 6)
(431, 430, 369, 60, 3) (433, 432, 379, 43, 5) (434, 433, 473, 38, 1) (435, 434, 466, 45, 3) (437, 436, 400, 34, 5)
(442, 441, 463, 20, 1) (442, 441, 469, 26, 2) (450, 447, 394, 16, 2) (450, 449, 487, 37, 2) (454, 453, 425, 28, 1)
(462, 459, 236, 34, 3) (462, 461, 484, 20, 1) (462, 461, 500, 37, 2) (463, 462, 496, 90, 4) (464, 463, 431, 105, 9)
(465, 463, 507, 55, 3) (471, 468, 536, 22, 2) (473, 471, 354, 57, 6) (479, 476, 426, 44, 6) (479, 478, 445, 56, 3)
(481, 476, 545, 24, 2) (482, 480, 437, 31, 3) (485, 482, 552, 22, 2) (487, 484, 423, 16, 3) (487, 486, 519, 35, 2)
(489, 487, 529, 19, 3) (492, 489, 560, 23, 1) (493, 491, 541, 24, 1) (494, 493, 522, 72, 4) (495, 494, 249, 115, 17)
(496, 493, 544, 16, 1) (497, 496, 465, 131, 9) (499, 497, 452, 60, 8) (500, 499, 533, 35, 2) (503, 501, 551, 20, 1)
(503, 501, 552, 25, 1) (509, 508, 477, 29, 2) (510, 507, 559, 17, 1) (511, 510, 540, 72, 5) (512, 511, 32, 32, 1)
(513, 511, 465, 59, 8) (513, 512, 545, 26, 3) (514, 511, 560, 28, 3) (516, 513, 572, 36, 3) (519, 518, 37, 46, 2)
(521, 518, 462, 42, 6) (521, 519, 474, 22, 2) (523, 522, 553, 28, 3) (524, 522, 487, 15, 2) (528, 525, 575, 27, 2)
(530, 529, 553, 22, 1) (532, 530, 133, 66, 1) (536, 532, 461, 30, 5) (536, 532, 593, 39, 4) (538, 536, 67, 67, 3)
(553, 551, 514, 16, 2) (555, 553, 499, 27, 2) (563, 560, 613, 26, 3) (565, 563, 516, 23, 2) (576, 574, 619, 21, 2)
(577, 575, 525, 33, 4) (581, 579, 535, 27, 3) (586, 583, 520, 45, 6) (586, 584, 73, 65, 3) (594, 592, 547, 28, 3)
(599, 597, 652, 30, 2) (608, 606, 521, 39, 3) (610, 607, 153, 48, 6) (627, 625, 573, 23, 3) (637, 635, 546, 40, 3)
(641, 639, 690, 21, 3) (654, 652, 595, 23, 4) (662, 660, 621, 17, 2) (668, 663, 739, 39, 4) (670, 667, 607, 20, 2)
(677, 674, 738, 55, 3) (681, 677, 742, 28, 2) (681, 678, 509, 53, 8) (687, 685, 735, 95, 5) (692, 690, -53, 27, 2)
(694, 692, 63, 63, 2) (695, 693, 652, 18, 2) (697, 686, 572, 18, 4) (697, 695, 746, 95, 5) (706, 704, 755, 48, 2)
(707, 705, 756, 47, 3) (719, 717, 652, 60, 7) (729, 726, 545, 56, 10) (730, 728, 781, 95, 5) (737, 734, 668, 62, 8)
(737, 735, 683, 25, 2) (740, 738, 687, 103, 10) (741, 738, 808, 54, 4) (742, 740, 695, 45, 3) (750, 747, 823, 24, 2)
(752, 749, 820, 54, 4) (752, 749, 844, 31, 2) (756, 754, 809, 95, 4) (758, 755, 851, 31, 1) (764, 762, 709, 101, 10)
(768, 765, 698, 21, 2) (777, 773, 703, 61, 7) (794, 789, 904, 22, 1) (795, 793, 737, 32, 3) (801, 796, 912, 22, 1)
(804, 802, 747, 27, 2) (808, 806, 749, 33, 3) (811, 808, 203, 66, 1) (816, 814, 771, 19, 2) (822, 817, 936, 23, 1)
(833, 830, 899, 73, 4) (839, 836, 902, 21, 1) (847, 845, 725, 81, 5) (850, 846, 773, 33, 3) (852, 849, 920, 71, 4)
(853, 851, 806, 20, 2) (855, 853, 794, 103, 9) (857, 854, 921, 22, 1) (884, 882, 821, 103, 9) (913, 911, 978, 90, 4)
(914, 911, 782, 45, 1) (921, 919, 858, 29, 2) (923, 920, 839, 32, 3) (931, 927, 696, 55, 3) (933, 930, 856, 39, 4)
(933, 931, 866, 105, 10) (943, 939, 1057, 58, 12) (943, 941, 876, 104, 10) (944, 942, 877, 56, 4) (978, 976, 911, 29, 2)
(986, 981, 875, 20, 4) (986, 984, 937, 21, 2) (993, 989, 900, 62, 8) (995, 991, 1085, 54, 4) (996, 992, 1093, 24, 1)
(1005, 1003, 1062, 72, 4) (1005, 1003, 1171, 81, 1) (1006, 1002, 1097, 54, 5) (1008, 1006, 943, 129, 11)
(1009, 1006, 918, 24, 4) (1009, 1007, 944, 131, 10) (1014, 1012, 1073, 29, 1) (1017, 1015, 948, 37, 3)
(1022, 1020, 957, 44, 3) (1027, 1025, 976, 22, 2) (1031, 1029, 1091, 29, 2) (1039, 1036, 1112, 47, 2)
(1039, 1037, 972, 127, 11) (1046, 1044, 975, 38, 3) (1064, 1061, 986, 25, 2) (1069, 1066, 1135, 41, 3)
(1072, 1069, 1146, 50, 2) (1076, 1073, 1164, 29, 2) (1088, 1085, 1164, 94, 6) (1089, 1086, 1166, 25, 1)
(1094, 1091, 1015, 102, 10) (1095, 1092, 1171, 48, 2) (1102, 1099, 1170, 40, 2) (1104, 1101, 1181, 97, 6)
(1109, 1105, 278, 67, 1) (1110, 1107, 1187, 48, 3) (1113, 1109, 279, 67, 1) (1136, 1133, 1054, 102, 10)
(1151, 1147, 1035, 27, 2) (1155, 1152, 1072, 26, 2) (1170, 1167, 1087, 52, 5) (1260, 1255, 1374, 55, 4)
(1273, 1270, 1351, 26, 1) (1282, 1279, 1373, 93, 4) (1293, 1290, 1372, 27, 1) (1319, 1316, 1225, 51, 5)
(1349, 1345, 1470, 30, 1) (1353, 1349, 1448, 46, 3) (1370, 1366, 1465, 97, 5) (1377, 1372, 1486, 36, 3)

| | | | |
|---|---|---|---|
| (1382, 1378, 1479, 95, 5) | (1389, 1384, 1499, 72, 5) | (1389, 1386, 1289, 108, 9) | (1393, 1389, 1489, 49, 4) |
| (1428, 1425, 1336, 130, 10) | (1460, 1457, 1360, 29, 2) | (1471, 1468, 1376, 128, 11) | (1472, 1469, 1594, 37, 2) |
| (1506, 1502, 1599, 41, 2) | (1514, 1508, 1651, 55, 4) | (1514, 1511, 1764, 81, 1) | (1613, 1609, 1512, 46, 3) |
| (1645, 1640, 1494, 59, 7) | (1647, 1643, 1529, 103, 10) | (1654, 1648, 1785, 72, 4) | (1667, 1661, 1799, 72, 5) |
| (1683, 1679, 1821, 69, 2) | (1690, 1685, 421, 83, 5) | (1694, 1689, 422, 84, 5) | (1695, 1691, 1834, 70, 2) |
| (1741, 1736, 1861, 49, 3) | (1760, 1755, 1882, 47, 2) | (1773, 1768, 1897, 96, 5) | (1812, 1808, 1941, 90, 4) |
| (1817, 1813, 1688, 104, 10) | (1865, 1861, 1633, 44, 5) | (1885, 1878, 235, 67, 3) | (1921, 1914, 2073, 72, 4) |
| (1932, 1921, 2142, 37, 2) | (1945, 1938, 2099, 72, 5) | (1975, 1971, 2087, 72, 4) | (1982, 1977, 1840, 103, 10) |
| (2173, 2166, 2346, 75, 5) | (2199, 2191, 2373, 72, 4) | (2248, 2242, 2405, 94, 6) | (2259, 2251, 2438, 72, 6) |
| (2290, 2281, 1715, 63, 3) | (2456, 2445, 1837, 53, 8) | (2457, 2451, 2281, 51, 5) | (2463, 2455, 2108, 44, 1) |
| (2487, 2479, 2685, 75, 5) | (2519, 2514, 2935, 81, 1) | (2560, 2555, 2709, 30, 1) | (2894, 2884, 3125, 71, 4) |
| (2911, 2903, 2701, 101, 10) | (2952, 2945, 3194, 70, 3) | (3067, 3059, 3318, 69, 3) | (3514, 3505, 3760, 96, 6) |
| (3602, 3593, 3344, 104, 9) | (3616, 3607, 3357, 104, 9) | (4148, 4139, 3107, 114, 5) | (4457, 4441, 4810, 73, 4) |
| (4644, 4631, 4966, 48, 2) | (10499,10471,11235,94,5) | (12774,12742,11859,104,9) | |